\documentclass{article}

\usepackage{arxiv}

\usepackage[utf8]{inputenc} % allow utf-8 input
\usepackage[T1]{fontenc}    % use 8-bit T1 fonts
\usepackage{hyperref}       % hyperlinks
\usepackage{url}            % simple URL typesetting
\usepackage{booktabs}       % professional-quality tables
\usepackage{amsfonts}       % blackboard math symbols
\usepackage{nicefrac}       % compact symbols for 1/2, etc.
\usepackage{microtype}      % microtypography

\usepackage{graphicx}

\usepackage{bm}
\newcommand{\myvec}[1]{\bm{#1}}

\newcommand{\revtwo}[1]{{#1}} 
\usepackage{amssymb} %% \mathbb
\usepackage{amsmath} %% \mathrm
\usepackage{multirow} %% \multirow command for table
\usepackage{booktabs} %% top/mid/bottom rule (in table)
\usepackage{siunitx}
\title{Calibration and Uncertainty Quantification of Convective Parameters in an Idealized GCM}

\author{Oliver R.~A.~Dunbar \\
Division of Geological and Planetary Sciences\\
California Institute of Technology\\
Pasadena, CA, USA\\
\texttt{odunbar@caltech.edu}
\And
Alfredo Garbuno-Inigo \\
Division of Actuarial Science, Statistics and Mathematics \\
Instituto Tecnol\'{o}gico Aut\'{o}nomo de M\'{e}xico, \\
Ciudad de M\'{e}xico, M\'{e}xico
\And
Tapio Schneider\\
Division of Geological and Planetary Sciences\\
California Institute of Technology\\
Pasadena, CA, USA\\
\And 
Andrew M.~Stuart\\
Division of Computational and Mathematical Sciences,\\
California Institute of Technology,\\
Pasadena, CA, USA\\
}

\begin{document}
\maketitle

\begin{abstract} % 250 words max; currently 243 words
Parameters in climate models are usually calibrated
manually, exploiting only small subsets of the available
data. This precludes both optimal calibration and quantification of uncertainties. Traditional Bayesian calibration methods that allow
uncertainty quantification are too expensive for climate
models; they are also not robust in the presence of internal climate variability. 
For example, Markov chain Monte Carlo (MCMC) methods typically
require $O(10^5)$ model runs and are sensitive to internal variability noise, rendering them infeasible for
climate models. Here we demonstrate an approach to model calibration and uncertainty
quantification that requires only $O(10^2)$ model runs and can accommodate internal climate variability. The approach consists of three stages:
(i) a calibration stage uses variants of ensemble Kalman inversion to calibrate a model by minimizing mismatches between model and data statistics; (ii) an
emulation stage emulates the parameter-to-data map with
Gaussian processes (GP), using the model runs in the
calibration stage for training; (iii) a sampling stage approximates the Bayesian posterior distributions by sampling the GP emulator with MCMC. We demonstrate the feasibility and
computational efficiency of this calibrate-emulate-sample (CES) approach in a perfect-model setting. Using an idealized general circulation model, we estimate parameters in a simple convection scheme from synthetic data generated with the model. The CES approach
generates probability distributions of the parameters that are good
approximations of the Bayesian posteriors, at a fraction of the
computational cost usually required to obtain them. Sampling from this 
approximate posterior allows the generation of climate predictions with quantified parametric uncertainties.
\end{abstract}

%%  --  --  --  --  --  --  --  --  --  --  --  --  --  --  --  --  --  --  --  --  --  --  --  --  %%
%
%  TEXT
%
%%  --  --  --  --  --  --  --  --  --  --  --  --  --  --  --  --  --  --  --  --  --  --  --  --  %%

%%% Suggested section heads:
% \section{Introduction}
%
% The main text should start with an introduction. Except for short
% manuscripts (such as comments and replies), the text should be divided
% into sections, each with its own heading.

% Headings should be sentence fragments and do not begin with a
% lowercase letter or number. Examples of good headings are:

% \section{Materials and Methods}
% Here is text on Materials and Methods.
%
% \subsection{A descriptive heading about methods}
% More about Methods.
%
% \section{Data} (Or section title might be a descriptive heading about data)
%
% \section{Results} (Or section title might be a descriptive heading about the
% results)
%
% \section{Conclusions}

\section{Introduction}
Calibrating climate models with available data and quantifying their uncertainties are essential to make climate predictions accurate and actionable. A primary source of uncertainties in climate models comes from representation of small-scale processes such as moist convection. Parameters in convection schemes and other parameterizations are usually calibrated by hand, using only a small fraction of data that are available. As a result, the calibration process may miss information about the small-scale processes in question. This paper presents a proof-of-concept, in an idealized setting, of how parameters in climate models can be calibrated using a substantial
fraction of the available data, and how uncertainties in the parameters can be quantified. We employ a new algorithm, called calibrate-emulate-sample (CES), which makes such calibration and uncertainty quantification feasible for computationally expensive climate models. CES reduces the hundreds of thousands of model runs usually required to quantify uncertainties in computer models to
hundreds, thereby achieving about a factor 1000 speedup. It leads to more robust calibration and uncertainty quantification in the presence of noise arising from chaotic variability of the climate system. We show how uncertainties in climate model parameters can be translated into quantified uncertainties of climate predictions through ensemble integrations.

%%  --  --  --  --  --  --  --  --  --  --  --  --  --  --  --  --  --  --  --  --  --  --  --  --  %%
%
%  TEXT
%
%%  --  --  --  --  --  --  --  --  --  --  --  --  --  --  --  --  --  --  --  --  --  --  --  --  %%

%%% Suggested section heads:
% \section{Introduction}
%
% The main text should start with an introduction. Except for short
% manuscripts (such as comments and replies), the text should be divided
% into sections, each with its own heading.

% Headings should be sentence fragments and do not begin with a
% lowercase letter or number. Examples of good headings are:

% \section{Materials and Methods}
% Here is text on Materials and Methods.
%
% \subsection{A descriptive heading about methods}
% More about Methods.
%
% \section{Data} (Or section title might be a descriptive heading about data)
%
% \section{Results} (Or section title might be a descriptive heading about the
% results)
%
% \section{Conclusions}

\section{Introduction}
%Type or paste your text here. The introduction gives a brief overview of the supporting information. You should include information %about as many of the following as possible (when appropriate):
% 1. a general overview of the kind of data files;
% 2. information about when and how the data were collected or created;
% 3. a general description of processing steps used;
% 4. any known imperfections or anomalies in the data.

The principal uncertainties in climate predictions arise from the
representation of unresolvable yet important small-scale processes,
such as those controlling cloud cover \cite{Cess89a, Cess90a, Bony05a,
  Stephens05, Bony06, Vial13a, Webb13b, Brient16b, Schneider17a}. 
These processes are represented by parameterization schemes,
which relate unresolved quantities such as cloud
statistics to variables resolved on the climate models' computational
grid, such as temperature and humidity. The parameterization
schemes depend on parameters that are a priori unknown, and so fixing the parameters is associated with uncertainty. The process of fixing these parameters to values that are most consistent with data is known as calibration, which generally requires solving an optimization problem.  Traditionally, however, parameters are calibrated (``tuned'') by hand, in a process that exploits only a small subset of the available observational data and relies on the knowledge and intuition of climate modelers about plausible ranges of parameters and their effect on the simulated climate of a model \cite{Randall97a, Mauritsen12a, Golaz13a, Hourdin13a, Flato13a, Hourdin17a,
  Schmidt17a, Zhao18b}. More recently, some broader-scale automated approaches that more systematically quantify the plausible range of parameters have begun to be explored \cite{Cou_etal20_pre,Hou_etal20_pre}. 
  
Opportunities to improve climate models lie in exploiting a larger
fraction of the available observational data together with high-resolution simulations, and learning from both, systematically and not manually \cite{Schneider17c}. To fully account for parametric uncertainty, we adopt a Bayesian view of the model-data relationship, which amounts to solving a Bayesian inverse problem. Given model, data, and prior information on the parameters, Bayesian inversion yields a posterior distribution of parameters. The mean or mode of the posterior distribution define the best parameter estimate, and the entire distribution provides uncertainty quantification through its spread about the mean or mode. We use the mismatch between climate statistics simulated with the model and those obtained from  observations or high-resolution simulations as the data likelihood in the Bayesian inverse problem to calibrate parameterizations in a climate model and to quantify their parametric uncertainties. Our focus is on learning from time-averaged climate statistics for three reasons: (1) time-averaged statistics are what is relevant for climate predictions; (2) time-averaged statistics vary more smoothly in space than atmospheric states, leading to a smoother optimization problem than that of atmospheric state estimation in numerical weather prediction (NWP); (3) time-averaging over long time-intervals reduces the effect of the unknown initial state of the system,
removing the need to determine it. Focusing on time-averaged climate statistics, rather than on instantaneous states or trajectories as in NWP, makes it possible to exploit climate observations and high-resolution simulations even when their native resolutions are very different from those of climate models.

While learning from climate statistics accumulated in time presents opportunities, it also comes with challenges. Accumulating statistics in time is computationally much more expensive than the forecasts over hours or days used in NWP. Therefore, we need algorithms for learning from data that minimize the number  of climate model runs required. Traditional methods for Bayesian calibration and uncertainty quantification such as Markov chain Monte Carlo (MCMC) typically require many iterations---often more than $10^5$---to reach statistical convergence e.g., \cite{book:Gey11}. Conducting so many computationally expensive climate model runs is not feasible, rendering MCMC impractical for climate model calibration \cite{Annan07a}. Additionally, while MCMC can be used to obtain the distribution of model parameters given data, it is not robust with respect to noise in the evaluation of the map from model parameters to data. Such noise, arising from natural variability in the chaotic climate system, can lead to trapping of the Markov chains in spurious, noise-induced local maxima of the likelihood function \cite{Cleary21a}. This presents additional challenges to using MCMC methods for climate model calibration. 

Here we demonstrate a new approach to climate model uncertainty quantification that overcomes the limitations of traditional Bayesian calibration methods in a relatively simple proof-of-concept. The approach---called calibrate-emulate-sample (CES) \cite{Cleary21a}---consists of three successive stages, which each exploit proven concepts and methods:
\begin{enumerate}
    \item In a calibration stage, we use variants of ensemble Kalman methods. These approaches were originally developed as a derivative-free method for state estimation \cite{evensen1994sequential,van1996data} and are now widely used in NWP \cite{Houtekamer16a}. The methods were subsequently developed for simultaneous
    state and parameter estimation, and variants were developed to deal with strongly nonlinear systems \cite{bocquet2012combining,gu2007iterative,li2007iterative,sakov2012,bocquet2014}. They were eventually recognized as a general purpose tool for the solution of inverse problems whose objective is parameter estimation \cite{oli,rey,reich2011dynamical,Evensen18a}. Here we use an optimization approach, referred to as ensemble Kalman inversion \cite{IglLawStu13}, which builds on the work of \cite{oli,rey}. \revtwo{Ensemble Kalman inversion may be understood as
    gradient descent projected to the ensemble subspace \cite{Schillings17a};
    relatedly, iterated EnKF smoothers may be interpreted as
    Gauss-Newton type methods confined to the ensemble subspace \cite{sakov2012}.}
    However, ensemble Kalman inversion and other ensemble Kalman methods do not provide a basis for systematic uncertainty quantification, except in linear Gaussian problems \cite{Annan07a,Gland09a,ErnSprSta15}. 
    
    \item In an emulation stage, we train an emulator on the climate model statistics generated during the calibration stage in order to quantify uncertainties. To emulate how the climate model statistics depend on parameters to be calibrated, we use Gaussian processes (GPs), a hierarchical machine learning method that learns smooth functions from a set of noisy training points \cite{Kennedy01a,Santner18a}. The GP approach also learns the statistics of the uncertainty in the resulting predicted function. The training points here are provided by the climate model runs performed in the calibration stage. 
    \item In a sampling stage, we approximate the
    posterior distribution on parameters given model and data,
    using the GP emulator to replace the parameter-to-climate statistics map. The emulator is used to estimate error statistics and to sample from the approximate posterior distribution with MCMC. Because the GP emulator is computationally cheap to evaluate and is smooth by virtue of the smoothing properties of GPs, this avoids the issues that limit the usability of MCMC for sampling from climate models directly. 
\end{enumerate}
The CES approach is described in detail in \cite{Cleary21a}, which provides a justification and contextualization of the approach in the literature on data assimilation and Bayesian calibration.  If uncertainty quantification is not required, the calibration step provides an effective, derivative-free parameter estimation method, and the emulation and sampling stages are not required. However, when uncertainty quantification is required, the role of the calibration stage is to provide a good set of training points in the vicinity of the posterior mean or mode, to train the emulator in the next stage of the algorithm.

The purpose of this paper is to demonstrate the feasibility of the approach for estimating parameters in an idealized general circulation model (GCM). This represents a proof-of-concept in a small parameter space and limited data space; how the methods scale up to larger problems will be discussed at the end.  

This paper is arranged as follows: Section \ref{sec:setup} describes the experimental setup, including the idealized GCM
and the generation of synthetic data from it. Section \ref{sec:methods} describes the CES approach
and the methods used in each stage. Section \ref{sec:results} describes the results of numerical experiments that use CES to calibrate parameters in the idealized GCM and quantify their uncertainties. It also demonstrates how sampling from the posterior distribution of parameters can be used to generate climate predictions with quantified uncertainties. Section \ref{sec:discussion} discusses and summarizes the results and their applicability to larger problems. 

\section{Experimental Setup}\label{sec:setup}

\subsection{General Circulation Model}

We use the idealized GCM described by \cite{Frierson06a} and \cite{OGorman08b}, which is based on the spectral dynamical core of the Flexible Modeling System developed at the Geophysical
Fluid Dynamics Laboratory. To approximate the solution of the hydrostatic primitive equations, it uses the spectral transform method in the horizontal, with spectral resolution T21 and with 32 latitude points on the transform grid. It uses finite differences with 10 unevenly spaced sigma levels in the vertical. We chose this relatively coarse resolution to keep our numerical 
experiments computationally efficient, so that comparison of CES with much more expensive methods is feasible. The lower boundary of the GCM is a homogeneous slab ocean (1~m mixed-layer thickness). Radiative transfer is represented by a semi-gray, two-stream radiative transfer scheme, in which the optical depth of longwave and shortwave absorbers is a prescribed function of latitude and pressure \cite{OGorman08b}, irrespective of the concentration of water vapor in the atmosphere (i.e., without an explicit representation of water vapor feedback). Insolation is constant and approximates Earth's annual mean insolation at the top of the atmosphere.

We focus our calibration and uncertainty quantification experiments on parameters in the GCM's convection scheme, which is a quasi-equilibrium moist convection scheme that can be viewed as a simplified version of the Betts-Miller convection scheme \cite{Betts86a,Betts86b,Betts93}. It relaxes temperature $T$ and specific humidity $q$ toward reference profiles on a timescale $\tau$ \cite{Frierson07b}:
\begin{equation}
    \frac{\partial T}{\partial t} + \cdots = - f_T\frac{T - T_\mathrm{ref}}{\tau}
\end{equation}
and
\begin{equation}
    \frac{\partial q}{\partial t} + \cdots = - f_T f_q \frac{q - q_\mathrm{ref}}{\tau}.
\end{equation}
Here, $f_T(z; T, q, p)$ is a function of altitude $z$ and of the thermodynamic state of an atmospheric column (dependent on temperature $T$, pressure $p$, and specific humidity $q$ in the column), which determines where and when the convection scheme is active; $f_q(T, q, p)$ is a function that modulates the relaxation of the specific humidity in non-precipitating (shallow) convection \cite{Frierson07b,OGorman08b}. The reference temperature profile is a moist adiabat, $T_\mathrm{ma}(z)$, shifted by a state-dependent and constant-with-height offset $\Delta T$, which is chosen to ensure conservation of enthalpy integrated over a column: $T_\mathrm{ref}(z) = T_\mathrm{ma}(z) + \Delta T$. The reference specific humidity $q_\mathrm{ref}(z)$ is the specific humidity corresponding to a fixed relative humidity RH relative to the moist adiabat $T_\mathrm{ma}(z)$. The two key parameters in this simple convection scheme thus are the timescale $\tau$ and the relative humidity \revtwo{parameter} RH. We demonstrate how we can learn about them from synthetic data generated with the GCM.

\subsection{Variable Selection and Generation of Synthetic Data}\label{sec:datmethod}

The idealized GCM with the simple quasi-equilibrium convection scheme has been used in numerous studies of large-scale atmosphere dynamics and mechanisms of climate changes, especially those involving the hydrologic cycle e.g.,  \cite{OGorman08b, OGorman08c, Bordoni08a, OGorman09a, Schneider10a, Merlis11a, OGorman11a, Kaspi11a, Kaspi13c, Levine15a, Bischoff14a, Wills17a, Wei18a}. We know from this body of work that the convection scheme primarily affects the atmospheric thermal stratification in the tropics, with weaker effects in the extratropics \cite{Schneider08c}. We also know that the relative humidity parameter (RH) in the moist convection scheme controls the humidity of the tropical free troposphere but likewise has a weaker effect on the humidity of the extratropical free troposphere \cite{OGorman11b}. Thus, we expect tropical circulation statistics to be especially informative about the parameters in the convection scheme. However, convection plays a central role in extreme precipitation events at all latitudes \cite{OGorman09a,OGorman09b}, so we expect statistics of precipitation extremes to be informative about convective parameters, and in particular to contain information about the relaxation timescale $\tau$.

\begin{figure}
\centering
  \noindent\includegraphics[width=0.5\textwidth]{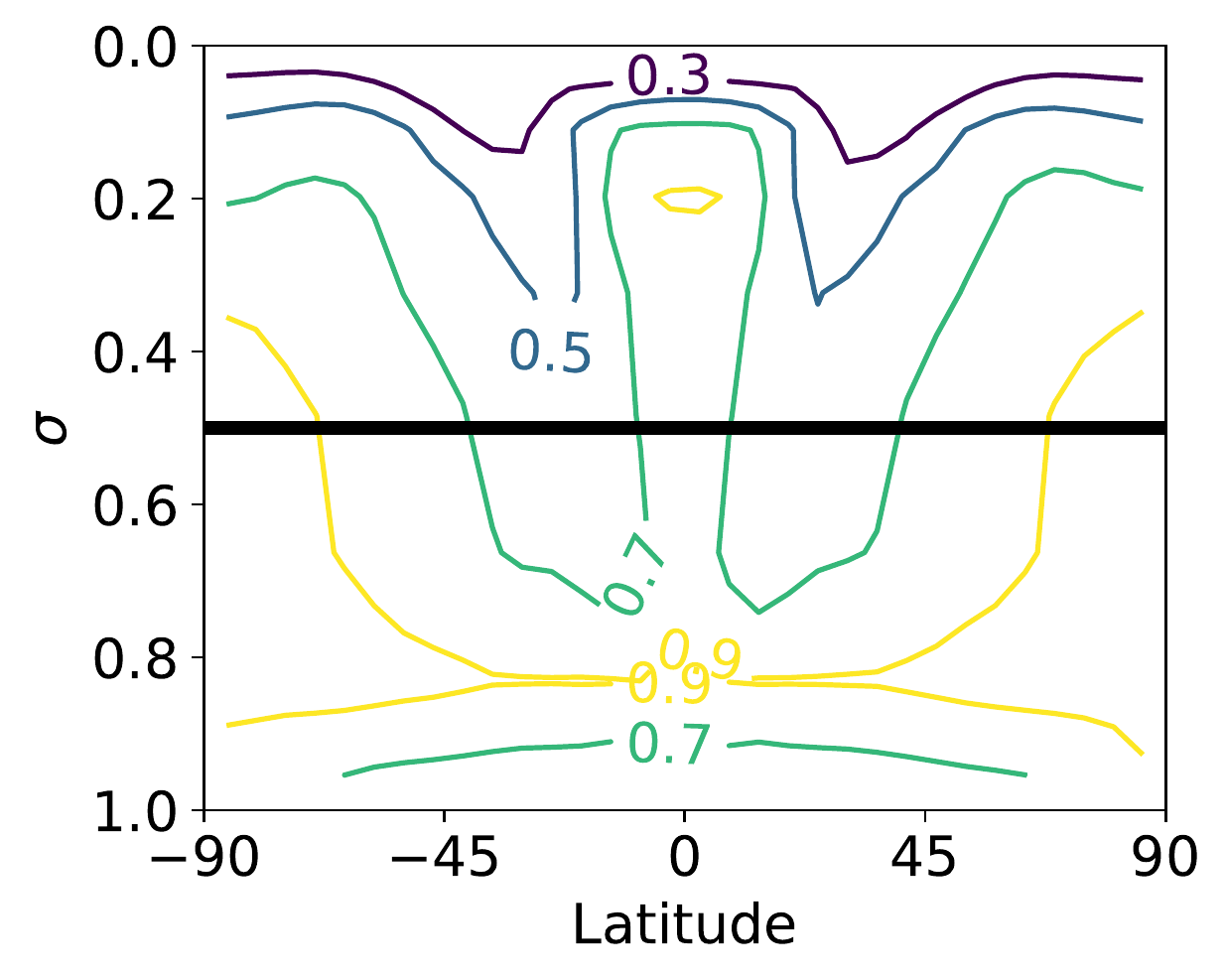}
  \caption{Zonal average of relative humidity averaged over one month. The black line shows the level at which data was extracted for computing data misfit functions.}
  \label{fig:rh_contours}
\end{figure}

\begin{figure}
  \includegraphics[width=\textwidth]{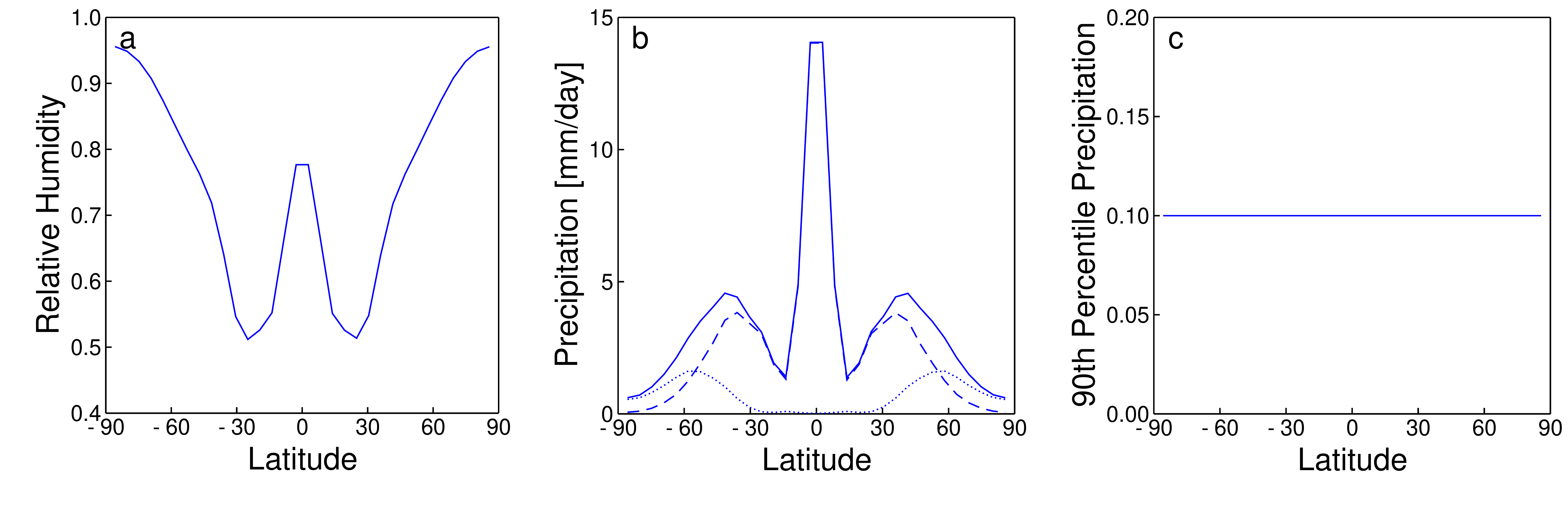}
  \caption{Long-term mean values of synthetic data. (a) Free-tropospheric relative humidity. (b) Total daily precipitation rate (solid) and its contributions from convection (dashed) and grid-scale condensation (dotted). (c) Probability of daily precipitation exceeding a 90th percentile (which is trivially 10\% in this case).}
  
  \label{fig:truth_means}
\end{figure}

As climate statistics from which we want to learn about the convective parameters, we choose 30-day averages of the free-tropospheric relative humidity, of the precipitation rate, and of a measure of the frequency of extreme precipitation. Because the GCM is statistically zonally symmetric, we take zonal averages in addition to the time averages. The relative humidity is evaluated at $\sigma=0.5$ (where $\sigma = p /p_s$ is pressure $p$ normalized by the local surface pressure $p_s$), as shown in Figure \ref{fig:rh_contours}. As a measure of the frequency of precipitation extremes, we use the probability that daily precipitation rates exceed a high, latitude-dependent threshold. The threshold is chosen as the latitude-dependent 90th percentile of daily precipitation in a long (18000 days) control simulation of the GCM in a statistically steady state. So for the parameters in the control simulation, the precipitation threshold is expected to be exceeded 10\% of the time at each latitude. The convective parameters in the control simulation are fixed at their reference values $\mathrm{RH}=0.7$ and $\tau = 2~\mathrm{h}$ \cite{OGorman08b}, and we collect the parameters in the vector $\myvec{\theta}^\dagger = (\theta_{\mathrm{RH}}^\dagger, \theta_\tau^\dagger)= (0.7, 2~\mathrm{h})$. Figure \ref{fig:truth_means} shows the mean relative humidity, the mean precipitation rate (broken down into its contributions coming from the convection scheme and from condensation at resolved scales), and the 90th percentile precipitation rate, from the control simulation averaged over 600 batches of 30-day windows. We use the single long control simulations of duration 18000 days only for the creation of Figure \ref{fig:truth_means} and for the estimation of noise covariances.

\subsection{Definition of noise covariance}
\label{sec:datcov}
\begin{figure}
\centering
  \noindent\includegraphics[height=6cm]{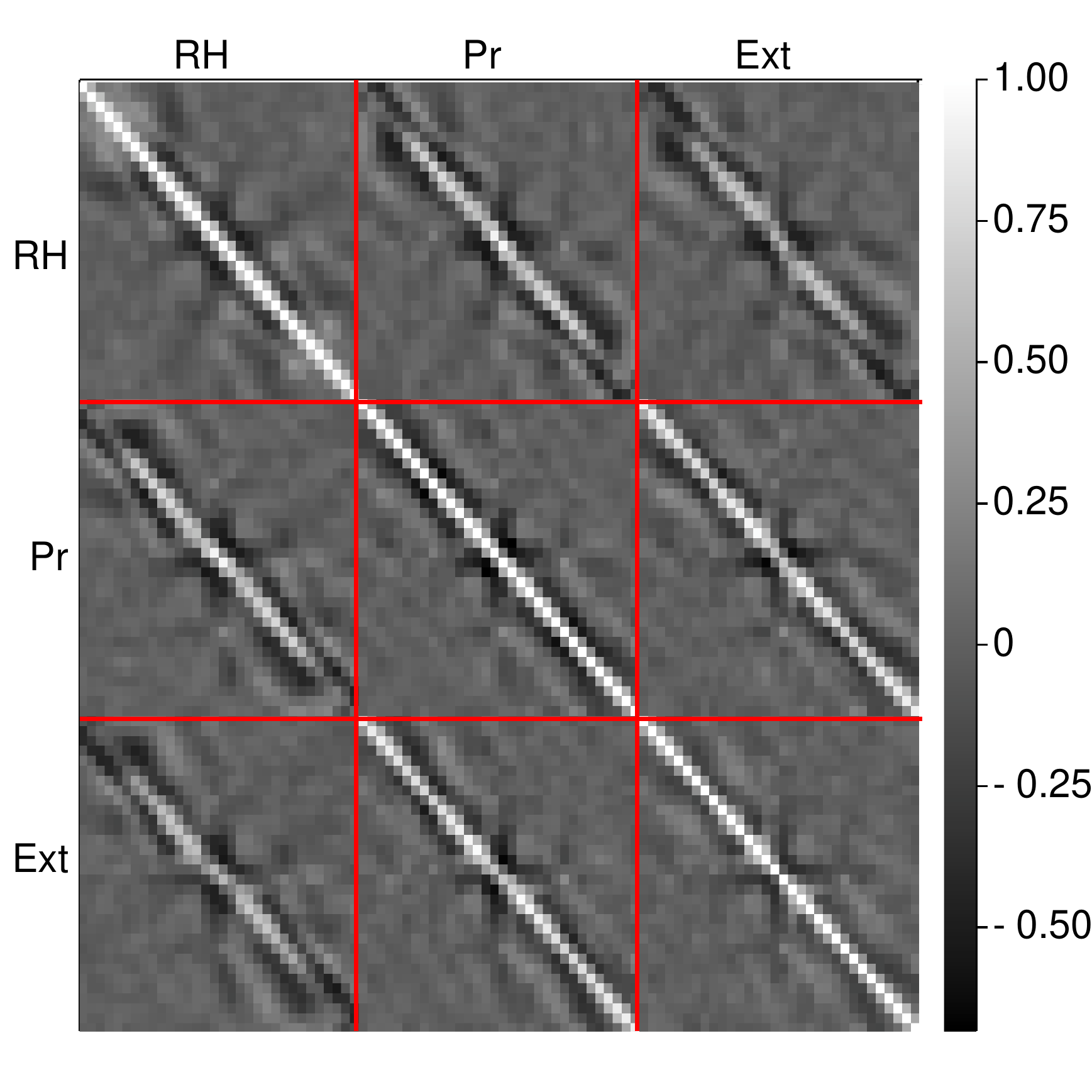}
  \caption{Correlation matrix associated with the internal variability estimated from the control simulation, plotted to illustrate the structure of $\Sigma$. The matrix blocks labeled (RH, Pr, Ext) are associated with the observed relative humidity, precipitation, and extreme precipitation.}
  \label{fig:sigma}
  \end{figure}
Estimation of model parameters requires specification of a noise covariance matrix, reflecting errors and uncertainties in the data. The principal source of noise in our perfect-model setting with synthetic data is sampling variability due to finite-time averaging with unknown initial conditions. The initial condition is forgotten at sufficiently long times because of the chaotic nature of atmospheric variability, so a central limit theorem quantifies the finite-time fluctuations around infinite-time averages that are caused by uncertain initial
conditions. Therefore, the asymptotic distribution of the fluctuations is a multivariate normal distribution $N(0,\Sigma(\myvec{\theta}))$ with zero mean and covariance matrix $\Sigma(\myvec{\theta})$. We estimate the covariance matrix at $\Sigma(\myvec{\theta}^\dagger)$, that is, with the parameters $\myvec{\theta}^\dagger$ in the control simulation. To estimate $\Sigma(\myvec{\theta}^\dagger)$, we run the GCM for 600 windows of length 30 days (because we use 30-day averages to estimate parameters) and calculate the sample covariance matrix of the 30-day means. With the 3 latitude-dependent fields evaluated at 32 latitude points, $\Sigma(\myvec{\theta}^\dagger)$ is a $96 \times 96$ symmetric matrix representing noise from internal variability in finite-time averages. Hereafter, we make the assumption that $\Sigma(\myvec{\theta}) \approx \Sigma(\myvec{\theta}^\dagger)$ for any $\myvec{\theta}$, and thus we treat $\Sigma$ as a constant matrix. In
practical implementations of this method, a corresponding constant $\Sigma$ can be estimated from climatology. We illustrate the correlation structure of $\Sigma$ in Figure~\ref{fig:sigma}.

\begin{figure}
  \includegraphics[width=\textwidth]{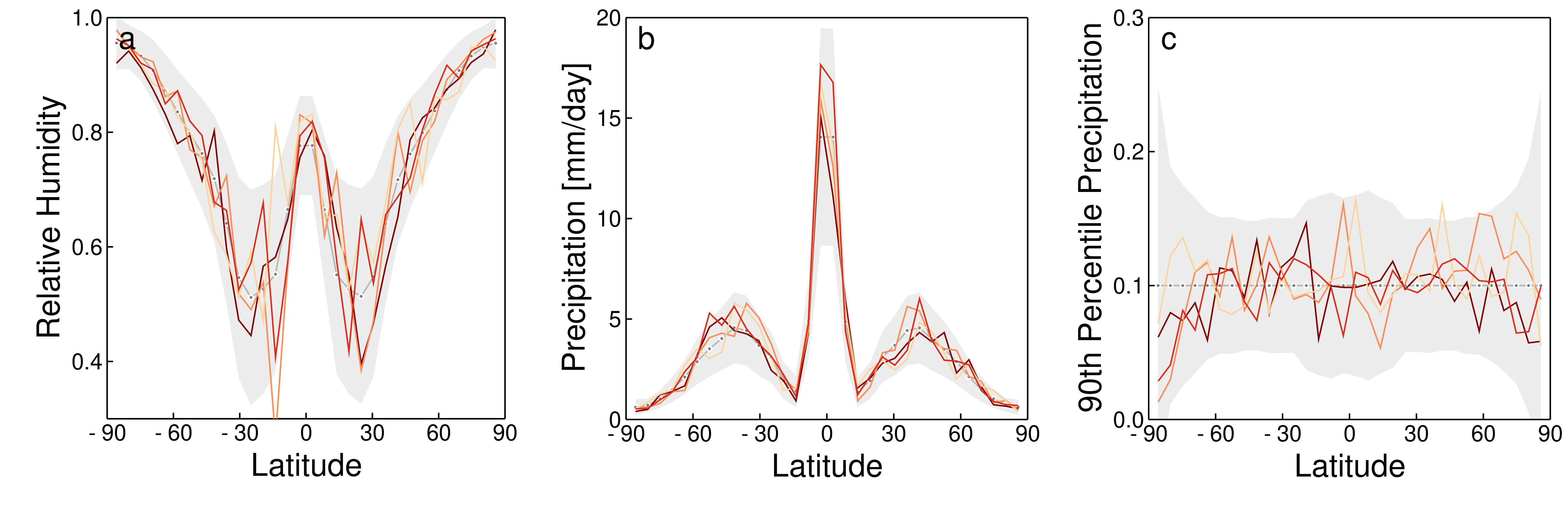}
  \caption{Four noisy realizations of the synthetic data, plotted in color over the underlying mean (grey circles) and 95\% confidence intervals from $\Gamma(\theta^\dagger)$ (grey bars). (a) Relative humidity. (b) Daily precipitation rate. (c) Probability of daily precipitation exceeding the 90th percentile of the long-term mean data.}
  \label{fig:truth_samples}
\end{figure}

To generate synthetic data, we also include the effect
of measurement error \cite{Kennedy01a}. We add Gaussian noise to the time-averaged model statistics, with a diagonal
covariance structure in data space.
We construct the measurement error covariance matrix $\Delta$ to be diagonal with entries $\delta_i>0$, where $i$ indexes over data type (the 3 observed
quantities) and latitude (32 locations). Combining this measurement covariance matrix $\Delta$ with the covariance matrix $\Sigma$ arising from internal variability leads to an inflated noise covariance matrix
\begin{equation}\label{eq:gammadef}
\Gamma = \Sigma+\mathop{\mathrm{diag}(\delta_i^2)} = \Sigma + \Delta.
\end{equation}
There are many options to pick $\delta_i$. We choose it by reducing the distance of the 95\% confidence interval to its nearest physical boundary for each $i$ by a constant factor $C$, so as to retain physical properties (e.g., precipitation must be nonnegative). Denote the mean $\mu_i$, variance $\Sigma_{ii}$, and a physical boundary set $\partial \Omega_i$ for each data $i$, we choose
\[
\delta_i = C \min\left(\text{dist}(\mu_i + 2\sqrt{\Sigma_{ii}}, \partial \Omega_i),\text{dist}(\mu_i - 2\sqrt{\Sigma_{ii}}, \partial \Omega_i)\right).
\]  
We take $C=0.2$. This value implies a significant noise inflation, with an average ratio of the standard deviations $\sqrt{\Gamma_{ii}}/\sqrt{\Sigma_{ii}}$ of $2.3$. Figure \ref{fig:truth_samples} shows the resulting data mean (grey circles), the 95\% confidence interval of the inflated covariance (grey ribbon), and four realizations of the data $\myvec{y}^{(1)},\dots ,\myvec{y}^{(4)}$ (yellow to red lines), each defined by taking a different 30-day average of the GCM, and adding a different realization of $N(0,\Delta)$. These four realizations will be used throughout when presenting our results. 

\section{Methods}\label{sec:methods}

\subsection{Misfit functions for time averaged data} \label{sec:objfunc}

Both calibration and uncertainty quantification in CES rely on a misfit function (standardized error) that quantifies mismatch between model output and data. Calibration minimizes a (possibly regularized) misfit function over the parameter space; uncertainty quantification samples from the posterior distribution, using a misfit function as the negative log-likelihood. To define the desired misfit function,
we introduce
$\mathcal{G}_T(\myvec{\theta};\myvec{z}^{(0)})$ and
$\mathcal{G}_\infty (\myvec{\theta})$, which
denote the mapping from the parameter vector $\myvec{\theta}$ to the $96$ data points, either
averaged over a finite time horizon ($T$) or over an
infinite time horizon ($\infty$). The former average
depends on the unknown initial condition $\myvec{z}^{(0)}$,
whereas the latter does not, because the initial condition is forgotten after a sufficiently long time.
We refer to $\mathcal{G}_T(\myvec{\theta};\myvec{z}^{(0)})$
as the forward model and $\mathcal{G}_\infty (\myvec{\theta})$ as the infinite time-horizon forward model.

To define the misfit function, we begin from the relationship between parameters $\myvec{\theta}$ and data $\myvec{y}$.
Expressed in terms of finite-time averages, this
relationship has the form
\begin{equation}\label{eq:invprob_inf0}
\myvec{y} = \mathcal{G}_T(\myvec{\theta};\myvec{z}^{(0)})+ \myvec{\eta},\qquad \myvec{\eta} \sim N(0,\Delta). 
\end{equation}
The data realizations $\myvec{y}^{(i)}$ for $i=1,\dots ,4$ in Section~\ref{sec:datcov} can therefore be seen as evaluations of \eqref{eq:invprob_inf0} for initial conditions $\myvec{z}^{(0)}_i$, noise realization $\myvec{\eta}_i\sim N(0,\Delta)$, and at a parameter $\myvec{\theta}^\dagger$. This form has the undesirable feature that it involves
$\myvec{z}^{(0)}$, a quantity that is not of intrinsic
interest. However, the central limit theorem asserts that 
\begin{equation} \revtwo{\label{eq:CLT}}
\mathcal{G}_T(\myvec{\theta};\myvec{z}^{(0)}) \approx \mathcal{G}_\infty(\myvec{\theta}) + \myvec{\sigma},\revtwo{\qquad \myvec{\sigma} \sim N(0,\Sigma).}
\end{equation}
This theorem quantifies the forgetting of the
initial condition after sufficiently long times, e.g., for the atmosphere, $T\gtrsim 15$ days \cite{BuiCheEmaMagSunZha19}. From the definition \eqref{eq:gammadef} of the total noise covariance matrix, $\Gamma = \Sigma + \Delta$, we can combine the observational and internal-variability noise and write
\begin{equation}\label{eq:invprob_inf}
\myvec{y} = \mathcal{G}_\infty (\myvec{\theta})+\myvec{\gamma}, \qquad \myvec{\gamma} \sim N(0,\Gamma). 
\end{equation}
We can equally interpret $\myvec{y}^{(i)}$ for $i=1,\dots ,4$ from Section~\ref{sec:datcov} as evaluations of \eqref{eq:invprob_inf} with noise $\myvec{\gamma}_i\sim N(0,\Gamma)$ and at parameter $\myvec{\theta}^\dagger$. This removes the dependence on initial condition but
is expressed in terms of infinite-time averages. Computing these averages directly is not
feasible, but (in the emulation phase) we introduce a procedure that enables us to learn a surrogate model for $\mathcal{G}_{\infty}$, using 
carefully chosen finite-time model evaluations $\mathcal{G}_T$.

 In the Bayesian approach to parameter learning, the aim
 is to determine the conditional distribution of parameters $\myvec{\theta}$ given a realization of data $\myvec{y}$ (written $\myvec{\theta}\mid \myvec{y}$)  by applying Bayes theorem, which states $\mathbb{P}(\myvec{\theta}\mid \myvec{y})$ is proportional to the product of the likelihood $\mathbb{P}(\myvec{y}\mid\myvec{\theta})$ and the prior $\mathbb{P}(\myvec{\theta})$. In the case where a surrogate
model for $\mathcal{G}_\infty$ is available, $\myvec{y}\mid\myvec{\theta}$ is defined as the pushforward of $\myvec{\theta}$ through a parameter-to-data map \eqref{eq:invprob_inf}, and we define the misfit function, the negative logarithm of the likelihood, as
\begin{equation}
  \Phi(\myvec{\theta}, \myvec{y}) = \frac{1}{2}\| \myvec{y} - \mathcal{G}_\infty(\myvec{\theta})\|_{\Gamma}^2\, ,
  \label{eq:objfunc_inf}
\end{equation}
where $\| \cdot \|_{\Gamma} = \| \Gamma^{-1/2} \cdot \|_2$ is the Mahalanobis distance. Before a surrogate
model for $\mathcal{G}_\infty$ is available, this function is infeasible to evaluate, but as $\myvec{y} \mid \myvec{\theta}$ can be defined using \eqref{eq:invprob_inf0} as well, we consider the related misfit function \revtwo{for the finite-time model evaluation at an initial condition $\myvec{w}^{(0)}$:}
\begin{equation}
  \Phi_T(\myvec{\theta}, \myvec{y};\revtwo{\myvec{w}^{(0)}}) = \frac{1}{2}\| \myvec{y} - \mathcal{G}_T(\myvec{\theta};\revtwo{\myvec{w}^{(0)}})\|_{\Gamma+\Sigma}^2\,.
  \label{eq:objfunc}
\end{equation}
\revtwo{In general, $\myvec{w}^{(0)}\neq \myvec{z}^{(0)}_i$, the initial conditions used to generate $\myvec{y}$}. We view evaluation of $\mathcal{G}_T$ from any initial condition as a random approximation of $\mathcal{G}_\infty$; hence, the additional internal-variability covariance matrix $\Sigma$ appears in \eqref{eq:objfunc}.
 
The CES algorithm in this context proceeds as follows: use optimization based 
on \eqref{eq:objfunc} to calibrate parameters; on the
basis of evaluations of $\mathcal{G}_T$ made during the
calibration, learn a
GP surrogate for $\mathcal{G}_\infty$; then 
use this surrogate to sample from the posterior
distribution of $(\myvec{\theta}\mid \myvec{y})$ defined using \eqref{eq:objfunc_inf}. We will henceforth neglect $\myvec{z}^{(0)}$ in our notation, and just write $\mathcal{G}_T(\myvec{\theta})$ and $\Phi_T(\myvec{\theta}, \myvec{y})$. Viewing the initial condition as random makes these objects  random as well. 

We have the following undesirable properties of the finite-time model average $\mathcal{G}_T(\myvec{\theta})$: (i) it is computationally expensive to evaluate for large $T$; (ii) it can be nondifferentiable or difficult to differentiate (e.g., because of non-differentiability of parameterization schemes in climate models); and (iii) evaluations of it are not deterministic (when one drops the explicit dependence on initial conditions). Our methodology, detailed in the upcoming sections, is constructed to overcome these difficulties.

An alternative approach to emulating the map $\mathcal{G}_\infty(\myvec{\theta})$ is to emulate \revtwo{the log-likelihood function} $\Phi(\myvec{\theta},\myvec{y})$ directly. This is often more computationally efficient as one always models a scalar function. \revtwo{But we found this efficiency comes at the cost of reduced accuracy and interpretability; for example, the choice of observational noise is unclear. In preliminary experiments, these drawbacks proved detrimental to performance, and so we abandoned the approach.}

 \subsection{Prior distributions}
 
 The priors of the physical parameters are taken to be the logit-normal and lognormal distributions, $\theta_{\mathrm{RH}}\sim \mathrm{Logitnormal}(0, 1)$ and
$\theta_{\tau}\sim \mathrm{Lognormal}(12~\mathrm{h},  (12~\mathrm{h})^2)$, for the relative humidity and timescale parameters, respectively. \revtwo{That is, we define the invertible transformation
\[
\mathcal{T}(\theta_\mathrm{RH}, \theta_\tau) = \left(\mathrm{logit}(\theta_\mathrm{RH}),~\ln\left(\frac{\theta_\tau}{1~\mathrm{s}}\right)\right),
\]
which transforms each parameter to values along the real axis.} Through an abuse of notation, we relabel the transformed (or computational) parameters as $\myvec{\theta}$, and the untransformed (or physical) parameters (relative humidity and timescale) are uniquely defined by $\mathcal{T}^{-1}(\myvec{\theta})$. The forward map $\mathcal{G}_T$ is defined as a composition $\mathcal{F}_T \circ \mathcal{T}^{-1}$ where $\mathcal{F}_T$ is the forward map defined on the physical parameters.

This choice allows us to apply our methods in the transformed space, where our priors are normally distributed and unbounded (namely logit $(\theta_{\mathrm{RH}}) \sim N(0,1)$ and log$(\theta_{\tau})\sim N(10.17,1)$); meanwhile the climate model is applied only to physically defined variables, $\theta_{\mathrm{RH}}\in[0,1]$ and $\theta_{\tau} \in [0,\infty)$. In this way the prior distributions enforce physical constraints on the parameters. 
 
\subsection{Calibrate: Ensemble Kalman Inversion}\label{sec:calibrate}

Ensemble Kalman inversion (EKI) \cite{IglLawStu13} is an
offline variant of ensemble Kalman filtering designed to
learn parameters in a general model, rather than states of
a dynamical system. EKI can be viewed as a derivative-free optimization algorithm. Given a set of data $\myvec{y}$, it iteratively evolves an ensemble of parameter estimates
both so that they achieve consensus and 
evolve toward the optimal parameter value $\myvec{\theta}^*$ (close to $\myvec{\theta}^\dagger$ if the quality/volume of the data
is good enough) that minimizes an objective function, in our case given by the misfit \eqref{eq:objfunc}, possibly with inclusion of a regularization term. It
has great potential for use with chaotic or stochastic models due to
its ensemble-based, derivative-free approach for optimizing
parameters. Theoretical work shows that noisy continuous-time versions of EKI exhibit an averaging effect that skips over 
fluctuations superimposed onto the ergodic averaged forward model \cite{duncan2021ensemble}. Empirical results suggest similar phenomena apply to EKI as implemented here, justifying the application to noisy forward model evaluations. Furthermore, the derivative-free approach scales well to high-dimensional parameter spaces, as evidenced by the use of Kalman filtering in numerical weather prediction, where billions of parameters characterizing atmospheric states are routinely estimated \cite{Kalnay03a}. This makes the
algorithm appealing for complex climate models. The algorithm is mathematically proven to find the optimizer,
within an initial, ensemble-dependent subspace,
for linear models \cite{SchStu17}. It is known to be effective for  high-dimensional nonlinear models \cite{IglLawStu13,jin1,jin2}, such as the nonlinear map from parameters to data represented by the idealized GCM we use in our proof-of-concept here.

The EKI algorithm we use is detailed in \cite{IglLawStu13}.
The algorithm iteratively updates an ensemble of parameters,
$\myvec{\theta}_m^{(n)}$, where $m=1,\dots M$ denotes an ensemble member, and the
superscript $n$ denotes the iteration count. The algorithm uses the ensemble to update
parameters according to the following equation
\[
  \myvec{\theta}_{m}^{(n+1)} = \myvec{\theta}_{m}^{(n)} + C_{\myvec{\theta} \mathcal{G}}^{(n)}\left( \Gamma + C_{\mathcal{G}\mathcal{G}}^{(n)}\right)^{-1}\left(\myvec{y} - \mathcal{G}_T(\myvec{\theta}_{m}^{(n)})\right)\,,
\]
where $C_{\mathcal{G}\mathcal{G}}$ is the empirical covariance of the
ensemble of quantities of interest from model runs, and $C_{\myvec{\theta}\mathcal{G}}$ is
the empirical cross-covariance of the ensemble of parameters and the ensemble of quantities of interest. The covariance matrix of the distribution of differences between realizations of $\myvec{y}$ and  $\mathcal{G}_T(\cdot)$ is $\Gamma+\Sigma$, which is approximated empirically in the update by $\Gamma+C_{\mathcal{G}\mathcal{G}}^{(n)}$.
When ensemble methods are used as approximate samplers \cite{oli,rey}, additional independent noise is added to $\myvec{y}$ at each iteration
and for every ensemble member; however, because
we are solving an optimization problem within this calibration phase, such noise is not added here.
%and for
%each ensemble member. However, because the individual %evaluations of $\mathcal{G}_T(\cdot)$ are affected by %internal variability, here we omit use of this additional %noise.

We initialize the algorithm by drawing an ensemble of size $M=100$ by sampling the parameter space from assumed prior
distributions on the parameters.

\subsection{Emulate: Gaussian Process Emulators (EKI-GP)}\label{sec:methods_gp}

\revtwo{In the emulation stage, we learn a surrogate for $\mathcal{G}_\infty$ by training a Gaussian process on samples of $\mathcal{G}_T$; the Gaussian process is a natural choice of tool here due to the Gaussianity of fluctuations in $\mathcal{G}_T$ about $\mathcal{G}_\infty$ described by the central limit theorem \eqref{eq:CLT}.} During the calibration stage with $N$ iterations and an ensemble of size $M$, we obtain a collection of input--output pairs
\[
\{\myvec{\theta}_m^{(n)}, \mathcal{G}_T(\myvec{\theta}_m^{(n)})\}, \qquad n=0,\dots, N, \quad m=1,\dots M.
\] 
The cloud of points $\{\myvec{\theta}_m^{(n)}\}$ from the calibration stage (a) spans the prior distribution, as initial EKI draws are from the prior, and (b) in later iterations, has a high density around the point $\myvec{\theta}^*$ to which EKI converges. We use regression to train a  GP emulator mapping $\myvec{\theta}$ to $\mathcal{G}_T(\myvec{\theta})$, using the input--output pairs $\{\myvec{\theta}_m^{(n)}, \mathcal{G}_T(\myvec{\theta}_m^{(n)})\}$, which are referred to as training points
in the context of GP regression. The emulation will be most accurate in regions with more training points, that is, around  $\myvec{\theta}^*$. This is typically near the true solution $\myvec{\theta}^\dagger$, and it is the region where the posterior parameter distribution will have high probability; this is precisely where uncertainty quantification requires accuracy. In effect, EKI serves as an effective algorithm for selecting good training points for GP regression. 

Gaussian processes emulate the statistics of the input--output pairs, using a Gaussian assumption. Specifically, we learn an approximation of the
form 
\[
  \mathcal{G}_T(\myvec{\theta}) \approx \mathcal{N}(\mathcal{G}_{\mathrm{GP}}(\myvec{\theta}), \Sigma_{\mathrm{GP}}(\myvec{\theta})).
\]
The approximation is learned from the input-output
pairs assuming that the outputs are produced from a mean function $\mathcal{G}_{\mathrm{GP}}(\myvec{\theta})$, and subject to normally distributed noise defined by a covariance function $\Sigma_{\mathrm{GP}}(\myvec{\theta})$, both dependent on the parameters.
The choice of notation here is to imply the fact that $\mathcal{G}_{\mathrm{GP}}(\myvec{\theta})$ serves to approximate the (unattainable) infinite-time average of the model $\mathcal{G}_{\infty}(\myvec{\theta})$, and $\Sigma_{\mathrm{GP}}(\myvec{\theta})$ serves to approximate the covariance matrix $\Sigma$.  Importantly, $\Sigma_{\mathrm{GP}}(\myvec{\theta})$ is $\myvec{\theta}$-dependent as it also includes the uncertainty in the approximation of the emulator at $\myvec{\theta}$ (for example, the emulator uncertainty $\Sigma_{\mathrm{GP}}(\myvec{\theta})$ will be large when $\myvec{\theta}$ is far from the inputs $\{ \myvec{\theta}_m\}$ used in training).  

The atmospheric quantities from which we learn about model parameters are correlated (e.g., relative humidity or daily precipitation at neighboring latitudes are correlated), resulting in a nondiagonal covariance matrix $\Sigma$. Any GP emulator therefore also requires a nondiagonal covariance $\Sigma_{\mathrm{GP}}(\myvec{\theta})$. We can enforce this, by mapping the correlated statistics from the GCM into a decorrelated space by using a  principal component analysis on $\Sigma$, and then training the GP with the decorrelated statistics to produce an emulator with diagonal covariance $\tilde{\Sigma}_{\mathrm{GP}}(\myvec{\theta})$. We use the notation $\tilde{(\cdot)}$ to denote variables in the uncorrelated space. To this end, we first decompose $\Sigma$ as
\[
  \Sigma = V D^2 V^T\,.
\]
Here, $V$ is an orthonormal matrix of eigenvectors of the covariance matrix $\Sigma$, and $D$ is the diagonal matrix of the square root of the eigenvalues, or the ordered standard deviations in the basis spanned by the eigenvectors of $\Sigma$. We store the outputs from the pairs as columns of a matrix $Y_{kl} = (\mathcal{G}_T(\myvec{\theta}_l))_k$, and then we change the basis of this matrix into the uncorrelated coordinates
\[
  \tilde{Y} = D^{-1}V^T Y\,.
\]
When trained on $\tilde{Y}$, GP returns the mean
$\tilde{\mathcal{G}}_{\mathrm{GP}}(\myvec{\theta})$ and the (diagonal)
covariance matrix $\tilde{\Sigma}_{\mathrm{GP}}(\myvec{\theta})$. We use tools from scikit-learn \cite{scikit-learn11} to train the emulator. After the diagonalization, we can train a scalar-valued GP for each of the 96 output dimensions, rather than having to train processes with vector-valued output. We construct a kernel by summing an Automatic Relevance Determination (ARD) radial basis function kernel and a white-noise kernel.  The ARD kernel is a standard squared exponential kernel, where each input dimension has an independent lengthscale hyperparameter.  This corresponds to regression,
rather than interpolation, and the variance of the
white noise kernel reflects the noise level assumed
in the regression. We train by learning 4 hyperparameters: the radial basis function variance, a lengthscale for each of the two parameters $\myvec{\theta}$ (due to ARD), and the white-noise variance. We train using the input--output pairs of the initial ensemble plus $N=5$ subsequent iterations of the EKI algorithm. We use $M=100$ ensemble members; thus, the training requires $(N+1)\times M = 600$ short (30-day) runs of our GCM. 

We continue using the uncorrelated basis in the sampling stage; where required, we transform the output of the emulator back into a correlated basis,
\begin{eqnarray}
  \mathcal{G}_{\mathrm{GP}}(\myvec{\theta}) &=& VD\tilde{\mathcal{G}}_{\mathrm{GP}}(\myvec{\theta}), \nonumber \\
  \Sigma_{\mathrm{GP}}(\myvec{\theta}) &=& VD\tilde{\Sigma}_{\mathrm{GP}}(\myvec{\theta})DV^T\,. \nonumber
\end{eqnarray}

\subsection{Sample: MCMC Sampling with a Gaussian Process Emulator}\label{sec:methods_mcmc}

To quantify uncertainties, we use MCMC to sample the
posterior distribution of parameters with the GP emulator.
The primary reason for using the GP emulator
goes back to the seminal paper by \cite{wynn} and concerns
the fact that it can be evaluated far more quickly than the GCM at a point in parameter space; this is important as we require more than $10^5$ samples within the likelihood $\mathbb{P}(\myvec{y} \mid \myvec{\theta})$ in a typical MCMC run to sample the posterior distribution of
parameters given data. However the emulator is
also important for two additional reasons:
(i) it naturally includes the approximation uncertainty (within $\tilde{\Sigma}_{\mathrm{GP}}$) of using an emulator; (ii) it smooths the likelihood function because we work with an approximation of \eqref{eq:objfunc_inf} based on the smooth $\mathcal{G}_{\infty}$, rather than \eqref{eq:objfunc} based on the noisy  $\mathcal{G}_T$; as a result, MCMC is less likely to get stuck in local extrema.

Recall that we trained the GP in uncorrelated coordinates. Within MCMC, one can either map back into the original coordinates or continue working in the uncorrelated space. We choose to continue working in the uncorrelated space, and so we map each data realization $\myvec{y}$ into this space: $\tilde{\myvec{y}} = D^{-1}V^T \myvec{y}$. In the Gaussian likelihood, we can use the GP emulated mean $\tilde{\mathcal{G}}_{\mathrm{GP}}(\myvec{\theta})$ and covariance matrix $\tilde{\Sigma}_{\mathrm{GP}}(\myvec{\theta})$ as  surrogates for the map $\mathcal{G}_\infty$ and the internal variability covariance  matrix $\Sigma$ (after passing to the uncorrelated coordinates). That is, we approximate the Bayesian posterior distribution as
\begin{align}
  \mathbb{P}(\myvec{\theta} \mid \tilde{\myvec{y}}) &\propto \mathbb{P}(\tilde{\myvec{y}}\mid\myvec{\theta}) \mathbb{P}(\myvec{\theta}) \nonumber\\
  &\propto \frac{1}{\sqrt{\det(\tilde{\Gamma}_{GP}(\myvec{\theta}))}}\exp\left(-\frac{1}{2}\| \tilde{\myvec{y}} - \tilde{\mathcal{G}}_{\mathrm{GP}}(\myvec{\theta})\|^2_{\tilde{\Gamma}_{\mathrm{GP}}(\myvec{\theta})}\right)\mathbb{P}(\myvec{\theta}) \nonumber\\
  &\propto \exp\left(-\frac{1}{2}\| \tilde{\myvec{y}} - \tilde{\mathcal{G}}_{\mathrm{GP}}(\myvec{\theta})\|^2_{\tilde{\Gamma}_{\mathrm{GP}}(\myvec{\theta})} - \frac{1}{2}\log \det \tilde{\Gamma}_{\text{GP}}(\myvec{\theta})  \right)\mathbb{P}(\myvec{\theta})\,. \label{eq:logposterior}
\end{align}
Here, $\tilde\Gamma_{\mathrm{GP}}(\myvec{\theta}) = \tilde \Sigma_{\mathrm{GP}}(\myvec{\theta})+ D^{-1}V^T\Delta V D^{-1} $ is the GP approximation of $\Gamma = \Sigma+\Delta$ in the uncorrelated coordinates. We include the (often overlooked) log-determinant term, arising from the normalization constant due to dependence of $\Gamma_{\text{GP}}$ on $\myvec{\theta}$. In the transformed parameter space, our prior $\mathbb{P}(\myvec{\theta})$ is also Gaussian and therefore can be factored inside this exponential, adding a quadratic penalty to the negative log-likelihood. The MCMC objective function is then defined to be the negative logarithm of the posterior distribution \eqref{eq:logposterior}; explicitly, given a Gaussian prior $N(\myvec{m},C)$ on the parameters, we define it as
\[
\Phi_{\mathrm{MCMC}}(\myvec{\theta},\tilde{\myvec{y}}) = \frac{1}{2}\| \tilde{\myvec{y}} - \tilde{\mathcal{G}}_{\mathrm{GP}}(\myvec{\theta})\|^2_{\tilde{\Gamma}_{\mathrm{GP}}(\myvec{\theta})} + \frac{1}{2}\log \det \tilde{\Gamma}_{\text{GP}}(\myvec{\theta})  + \frac{1}{2}\|\myvec{\theta}- \myvec{m}\|^2_{C}.
\]
Evaluation of $\Phi_{\mathrm{MCMC}}$ requires only the
mean and covariance from the GP emulator. Furthermore $\Phi_{\mathrm{MCMC}}$ is smooth and so is suitable for use within an MCMC algorithm to generate samples from the approximate posterior distribution of the parameters. \cite{Cleary21a} contains further discussion of MCMC using GPs to emulate
the forward model, including situations where
data comes from finite time-averages but the emulator
is designed to approximate the infinite time-horizon forward model.

We use the random walk metropolis algorithm for MCMC sampling. The
priors chosen were the same, physics-informed priors used to
initialize EKI. We choose the proposal distribution also as a Gaussian with covariance proportional to the prior covariance. 
%The proposal covariance \hl{[might be helpful to explain what this covariance is. There are a lot of covariance matrices here, and it's not easy for readers to keep them straight]} was initially estimated with the empirical covariance of iteration 5 from EKI, and updated during a 10,000 sample burnin phase. 
The MCMC run consists of a burn-in of 10,000 samples followed by 190,000
samples. 
%In our experiments, we ran MCMC individually for each of the
%300 samples \hl{[why 300? explain where they come from]} of the synthetic truth.

\subsection{Benchmark Gaussian process (B-GP)}
The performance of any emulator is dependent on the training points. Since we use an adaptive procedure (EKI) to concentrate the training points, which is the novel approach
introduced in \cite{Cleary21a}, we also train a benchmark emulator to compare our results with those resulting from more traditional, brute-force approaches to the emulation. 

For this purpose, we use a GP emulator trained on a uniform set of points. Even in two dimensions, it is prohibitively costly for this set to span the support of the prior distribution. Instead, we use knowledge of the location and size of the posterior distribution to place a uniform grid of $40 \times 40 = 1600$ training points over $[-1.25,-0.5]\times [8.0,10.0]$ in the transformed parameter space. This corresponds to  $[0.62,0.77]\times [0.83~\text{h}, 6.12~\text{h}]$ in the untransformed parameter space and captures the region of high probability of the posterior. The use of posterior knowledge here reduces the number of training points by a factor of 20 when compared to spanning 95\% of the prior mass; no such posterior information is used in the CES algorithm, which automatically places training points where they are needed. The benchmark emulator uses the same kernel and training setup as in section \ref{sec:methods_gp}, and we use the trained emulator in
MCMC experiments in the same way as described in Section \ref{sec:methods_mcmc}. To distinguish the two methods, we refer to the EKI-trained GP as EKI-GP and the benchmark (traditionally trained) GP as B-GP. 

\section{Results} \label{sec:results}

To demonstrate the dependence of the parameter uncertainty on the realization of the (inflated) synthetic data, we reproduce the experiments 4 times with each of the four realizations shown in Figure \ref{fig:truth_samples}. We denote these four sets of data $\myvec{y}^{(1)},\dots, \myvec{y}^{(4)}$.

\subsection{Calibrate: Ensemble Kalman Inversion}
\begin{figure}
\centering
  \noindent\includegraphics[width=0.65\textwidth]{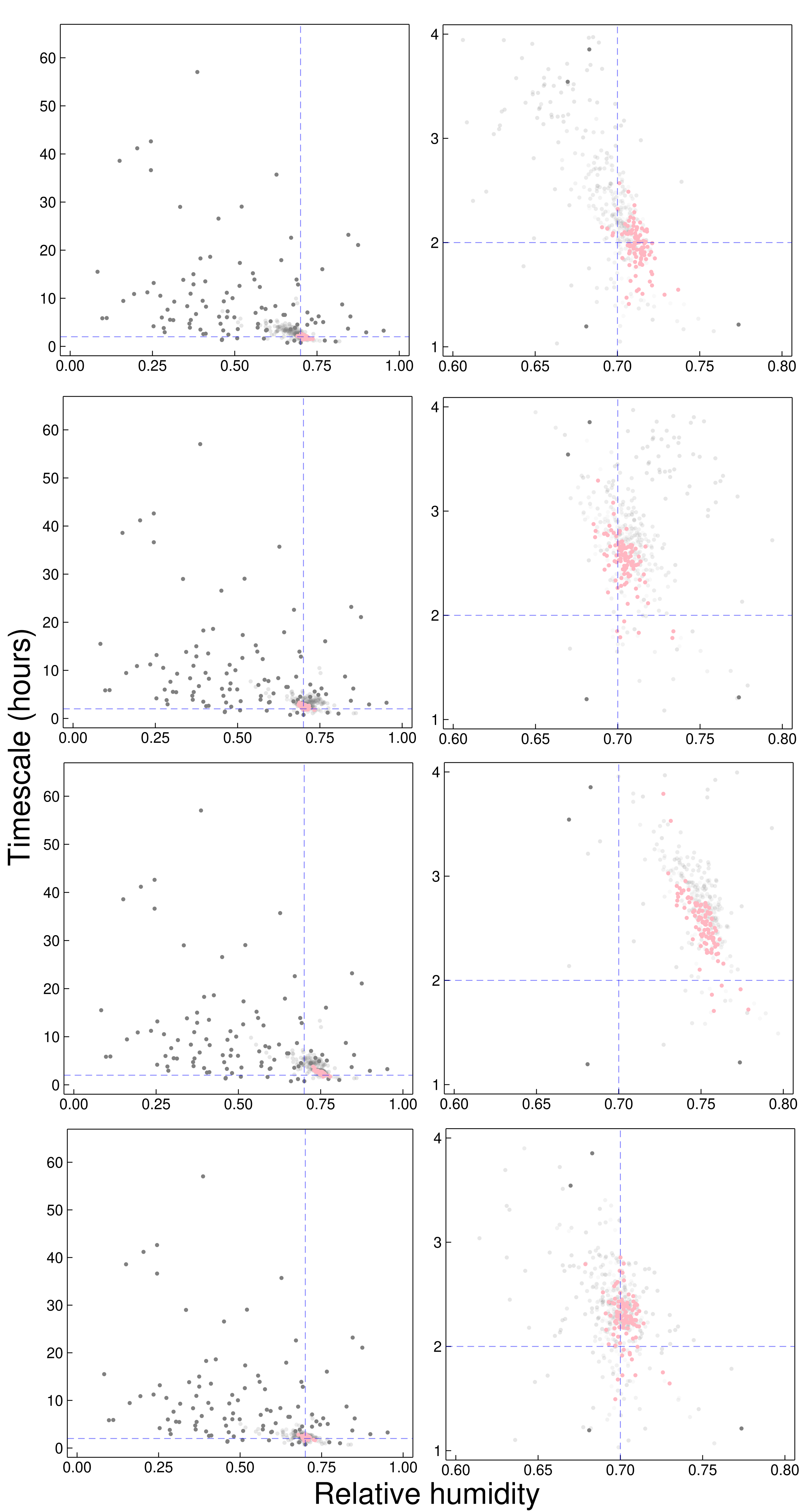}
  \caption{EKI ensemble at iterations 0 to 5 displayed as particles in parameter space. Left column:  all members; right column: zoom-in near true parameter values. Each row represents optimization with a different data vector $\myvec{y}^{(i)}$ from Figure~\ref{fig:truth_samples}. The (initial) prior ensemble 0 is highlighted in dark grey, and the final ensemble 5 is highlighted in pink. The intersection of the dashed blue lines represents the true parameter values used to generate observational data from the GCM.}
  \label{fig:eki}
\end{figure}

\begin{figure}
  %% \setfigurenum{S3}
  \centering
  \noindent\includegraphics[width=0.8\textwidth]{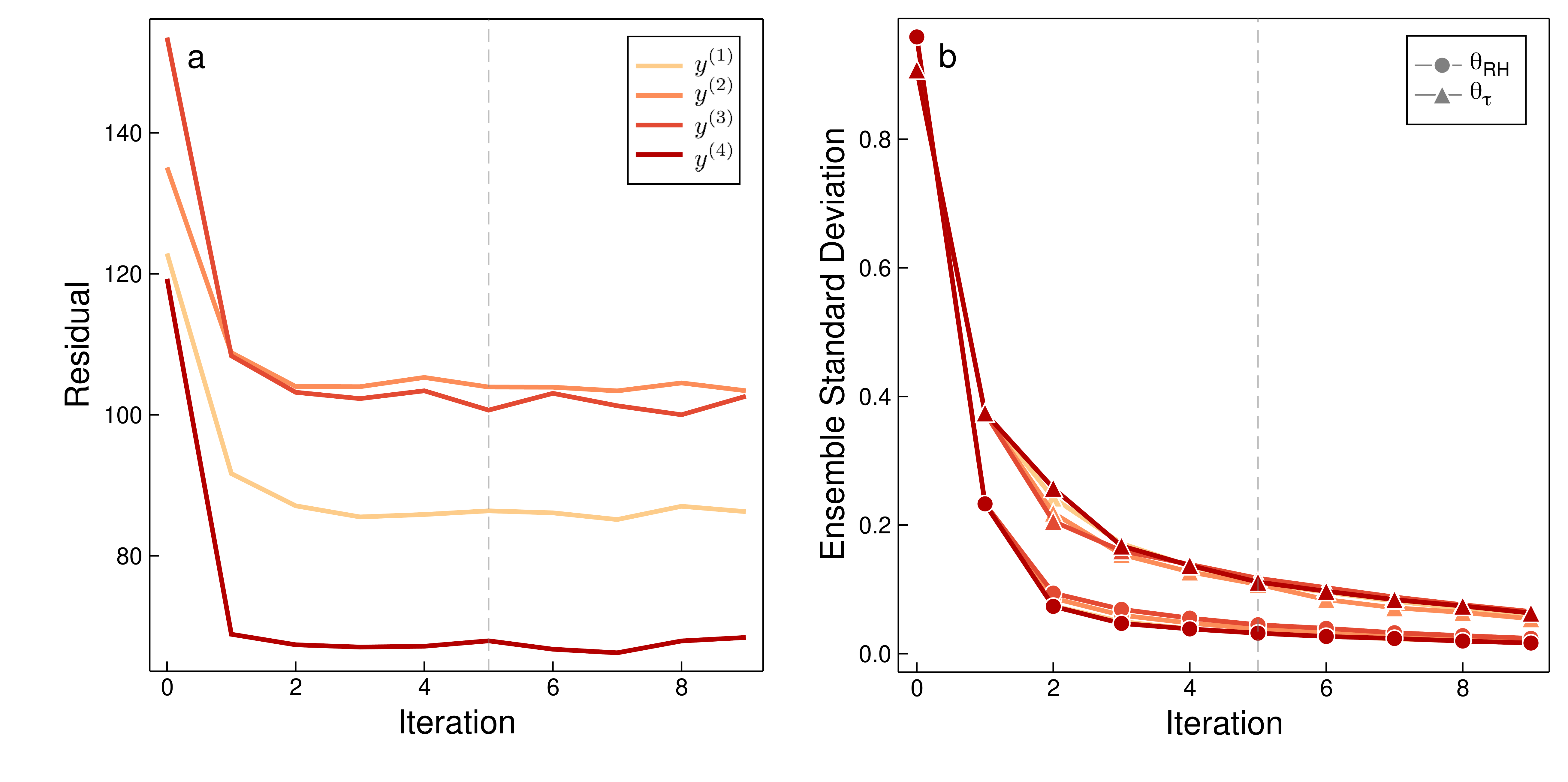}
  \caption{Convergence behaviour tests over 9 iterations of EKI for each realization of the data. The vertical dashed line marks the final iteration of Figure \ref{fig:eki}; we also show behaviour of 4 further iterations. (a) Ensemble-mean residuals relative to synthetic data for each EKI iteration. (b) Standard deviation of ensemble for the relative humidity parameter (circle) and timescale parameter (triangle) for each realization.}
  \label{fig:resid_stddev}
\end{figure}
We use the first 6 iterations of EKI in the CES training process. The initial
ensemble is spread over the whole parameter space but collapses within a few iterations near the true parameter values---to within 10\% error in $\theta_{\mathrm{RH}}$ and 30 minutes error in $\theta_{\tau}$ (Figure \ref{fig:eki}). That is, the algorithm evolves toward
consensus and toward the true solution. Biases arise from the realization of internal variability and the
realization of the observational noise in each $\myvec{y}^{(\cdot)}$.

To check for EKI convergence, we evaluate an additional 4 EKI iterations (labeled 0 to 9). At each iteration $n$, we compute residuals of the ensemble mean for each
realization of the synthetic data $\myvec{y}^{(1)},\dots,\myvec{y}^{(4)}$ created at the true parameters $\myvec{\theta}^\dagger$,
\begin{equation*}
  \mathrm{Residual}(n;\myvec{y}^{(i)}) = \left\lVert \frac{1}{M}\sum_{m=1}^M \mathcal{G}_{T}(\myvec{\theta}^{(n)}_m) - \myvec{y}^{(i)} \right\rVert_{\Gamma}^2\,,
\end{equation*}
weighting the residuals by the covariance matrix $\Gamma$ of the synthetic data. Figure
\ref{fig:resid_stddev}(a) shows the residual as a function of EKI iteration. The residual decreases quickly over the first few iterations, before plateauing for subsequent
iterations. Figure
\ref{fig:resid_stddev}(b) shows standard deviations of the ensemble of parameters. The standard deviations decrease monotonically from iteration to iteration, reflecting the evolution
toward consensus. The behavior is qualitatively similar for all realizations; quantitative differences reflect different realizations of internal variability in the different data realizations. This behavior reflects the fact that EKI is
an optimization method for calibrating parameters: it is not
constructed to learn uncertainty but rather to reach consensus around a single parameter value that makes the misfit small
\cite{Schillings17a}.

\subsection{Emulate: Validation}

\begin{figure}
  %% \setfigurenum{S3}
  \centering 
  \includegraphics[width=0.65\textwidth]{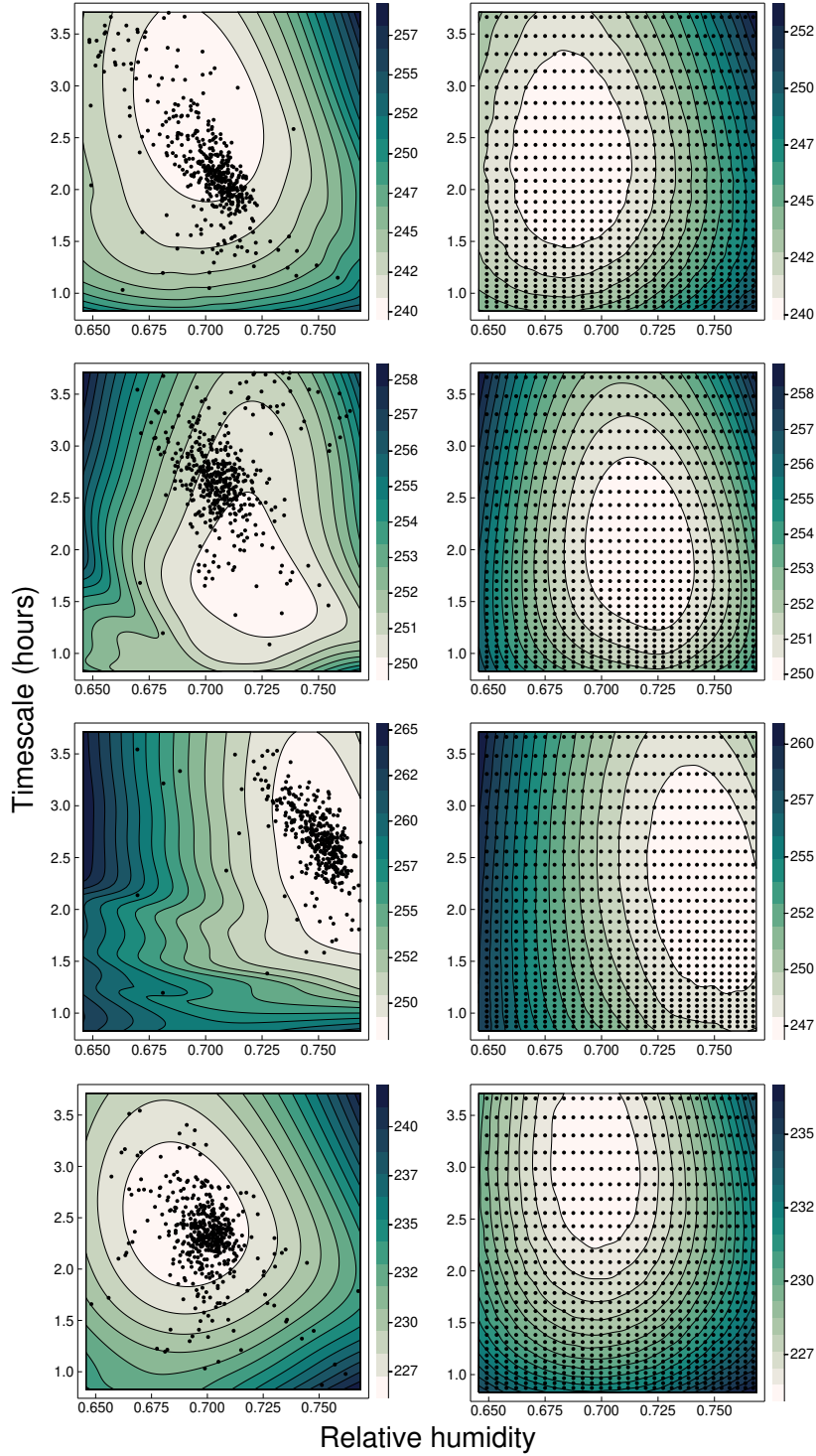}
  \caption{Training points for the GP emulators (EKI-GP and B-GP), plotted over the objective function $\Phi_{\mathrm{MCMC}}$ for different data realizations $\myvec{y}^{(1)},\dots,\myvec{y}^{(4)}$ (rows). Left column: particles representing members of the first 6 EKI iterations. Right column: grid (uniform in the transformed parameters) used to train the benchmark Gaussian process. In both cases, some additional training points fall outside of the plotting domain. }
  \label{fig:tp}
\end{figure}

\begin{figure}

  \centering
  \includegraphics[width=0.9\textwidth]{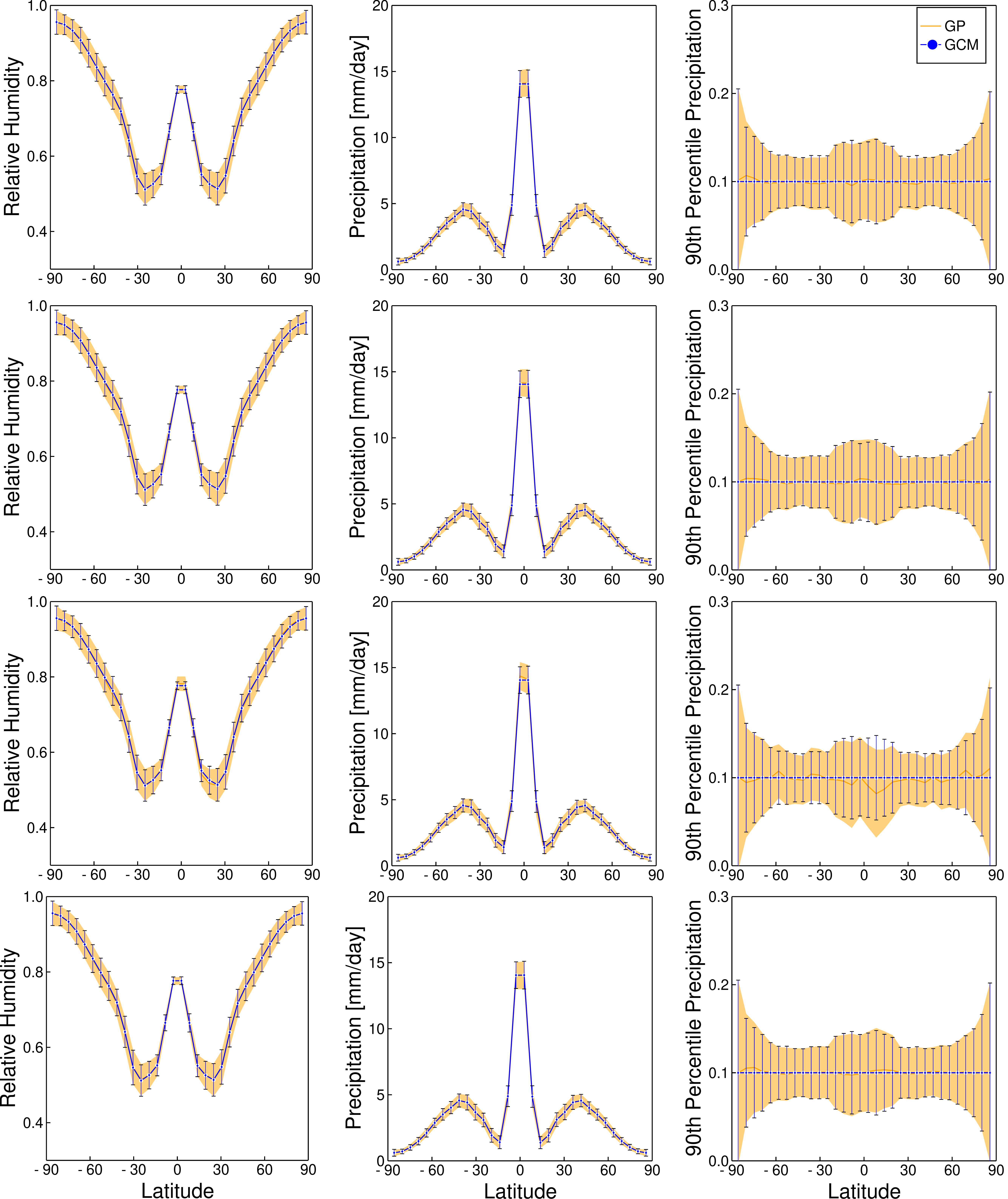}
  \caption{Comparison between the GCM statistics at the true parameters $\myvec{\theta}^\dagger$ and the trained EKI-GP emulator at $\myvec{\theta}^\dagger$. The four rows correspond to using EKI with the data vectors  $\myvec{y}^{(1)},\dots,\myvec{y}^{(4)}$. Blue lines: GCM mean (dots) averaged over 600 30-day runs, with the error bars marking a 95\% confidence interval from variances on the diagonal of $\Gamma$. Orange: predicted mean (line) and 95\% confidence interval (shaded region) produced by the GP emulator. }
  \label{fig:validation}
\end{figure}

\begin{figure}

\centering
\includegraphics[width=0.9\textwidth]{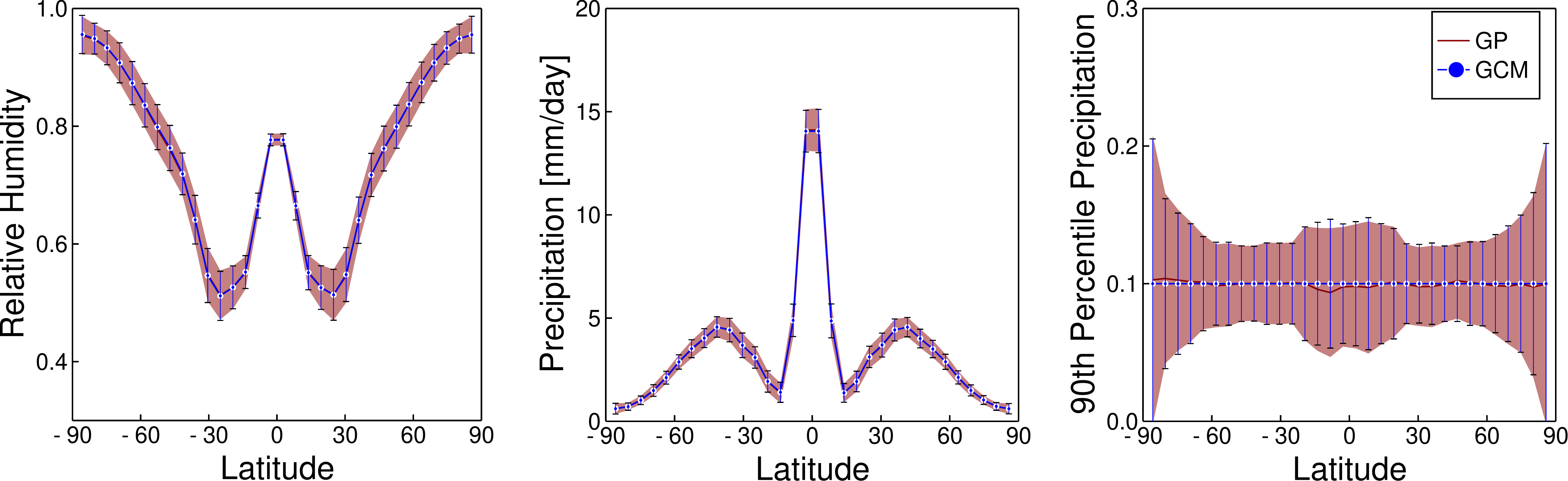}

\caption{Comparison between the GCM statistics at the true parameters $\myvec{\theta}^\dagger$ and the trained B-GP emulator predictions at $\myvec{\theta}^\dagger$. Blue: GCM mean (dots) averaged over 600 30-day runs, with the error bars marking a 95\% confidence interval from variances on the diagonal of $\Gamma$. Dark red: predicted mean (line) and 95\% confidence interval (shaded region) produced by the B-GP emulator.  }
  \label{fig:validation_benchmark}
\end{figure}
Figure \ref{fig:tp} shows the parameter values used for training points for the EKI-GP and B-GP. We use the first 6 EKI iterations (i.e., 600 training points) for training. These are plotted over the associated objective function $\Phi_{\mathrm{MCMC}}$. The panels in the left column correspond to the EKI-GP using data vectors $\myvec{y}^{(1)},\dots, \myvec{y}^{(4)}$. We see that the resulting objective  functions  $\Phi_{\mathrm{MCMC}}$ from training EKI-GP at the marked training points (black dots) lead to unimodal distributions with a minimum near to a significant number of the training points; there are also training points that fall outside of the plotting domain (see Figure \ref{fig:eki} for their extent). The right column of Figure \ref{fig:tp} shows the benchmark grid for B-GP, which we use as a comparison for the EKI-GP method; the contours of $\Phi_{\mathrm{MCMC}}$ were calculated using the same realization as their counterpart EKI-GPs. We see that for each realization the EKI-GP and B-GP produce objective functions that are qualitatively similar in terms of the magnitude of the minimum, the location of the minimum, and the approximate shape of the objective function; the quantitative differences are accounted for by differing geometry and density of training points (and hence a difference in approximation uncertainty). In both settings, the objective function $\Phi_{\mathrm{MCMC}}$ is smooth because the GP smoothly approximates $\mathcal{G}_{\infty}$. 

EKI-GP shows similar results for the objective function as B-GP, at a fraction of the computational effort. B-GP is far less practical as
a methodology than is EKI-GP because it does not scale
well to high-dimensional parameter spaces; it requires  many more
training points than EKI-GP. Even in our two-dimensional experiments, we needed to use posterior information to reduce the number of training points for B-GP by a factor of 20. The B-GP comparison is included simply
to demonstrate that EKI-GP achieves comparable results to
those achieved by means of traditional emulation.

We validate the emulator approximation to the data by making a prediction at the true parameters $\myvec{\theta}^\dagger$. We display $\mathcal{G}_{\mathrm{GP}}(\myvec{\theta}^\dagger)$ and the 95\% confidence intervals computed using the variance from $\Sigma_{\mathrm{GP}}(\myvec{\theta}^\dagger)$ in Figure \ref{fig:validation} for EKI-GP, and in Figure \ref{fig:validation_benchmark} for B-GP. The rows of Figure \ref{fig:validation} correspond to the EKI-GP results for $\myvec{y}^{(1)},\dots, \myvec{y}^{(4)}$. In both figures we also show the statistics of 600 30-day samples from the control simulation at $\myvec{\theta}^\dagger$.  Both the mean and 95\% confidence intervals of all EKI-GP emulators (orange line and ribbon) closely match the statistics from the GCM runs (blue dots and error bars), as does the prediction from the B-GP (dark red line and ribbon). The training data are sufficient to ensure that the predicted 95\% confidence interval from the emulators do not produce unphysical values (such as giving negative precipitation rates, or relative humidities outside $[0,1]$).

\subsection{Sample: MCMC Sampling}
\begin{figure}
\centering
  \includegraphics[width=0.65\textwidth]{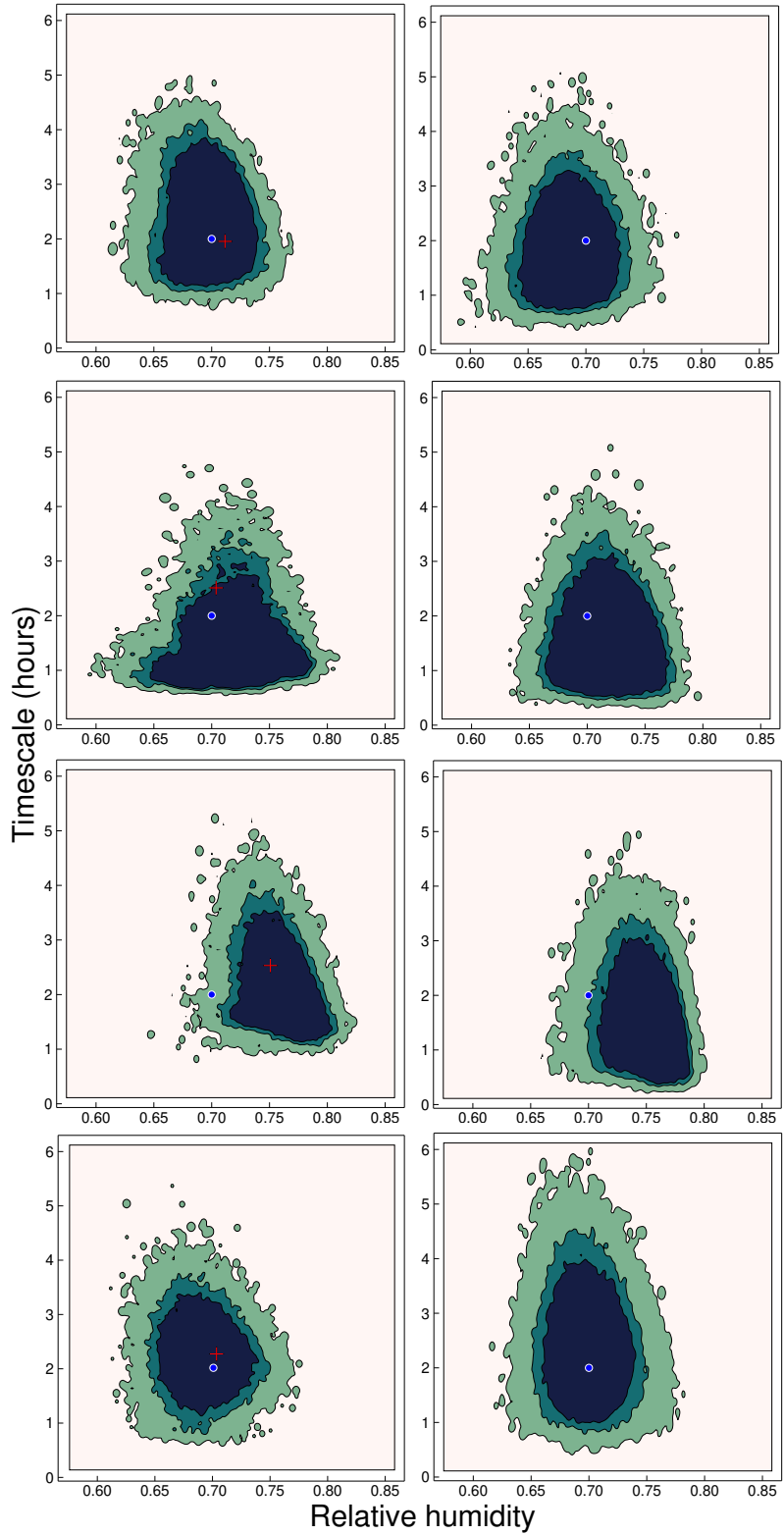}
  \caption{Density plot of MCMC samples of the posterior distribution. The contours are drawn to contain 50\%, 75\%, and 99\% of the distribution generated from the samples. The left column show distributions learned using EKI-GP at  $\myvec{y}^{(1)},\dots,\myvec{y}^{(4)}$, and the right column using B-GP at the same data realizations. The blue dot represents the true parameters, while the red plus is the average across particles in the 6th EKI iteration. }
  \label{fig:mcmc}
\end{figure}

\renewcommand{\arraystretch}{1.5}
\begin{table}[]
\centering
%%%%%%% THIS TABLE CONTAINS THE SAMPLES
%   \begin{tabular}{lccc}
%     \toprule
%     \multirow{2}{*}{} & {$\sigma_{RH}$} & {$\sigma_{\tau}$ (hrs)}\\
%     \midrule
%     EKI (Iteration 9)  & 0.014886, 0.0194042, 0.0207519, 0.0142504   & 0.0534855, 0.0490314, 0.057362, 0.0551679

% \\
%     MCMC (EKI-GP)   &     0.0888076, 0.123785, 0.0965848,   0.0890853 & 0.272798, 0.3077833, 0.244517, 0.233396\\
%     MCMC (B-GP) & 0.091078, 0.102536, 0.0967786, 0.0924501& 0.322513, 0.362513, 0.431368, 0.321089\\
%     \bottomrule\\
%   \end{tabular}
%%%%%%%% THIS TABLE CONTAINS THE AVERAGES
    \begin{tabular}{lccc}
    \toprule
    \multirow{2}{*}{} & {$\sigma_{RH}$} & {$\sigma_{\tau}$ (hrs)}\\
    \midrule
    EKI (Iteration 9)  & 0.017  & 0.053

\\
    MCMC (EKI-GP)   &     0.099 & 0.265\\
    MCMC (B-GP) & 0.096 & 0.359\\
    \bottomrule\\
  \end{tabular}
  \caption{Average standard deviations of parameters from EKI and MCMC experiments over $\myvec{y}^{(1)},\dots,\myvec{y}^{(4)}$. }
  \label{table:stddev}
\end{table}

We use an MCMC algorithm to generate a set of samples from the posterior distribution with the help of the GP emulator. We choose the random walk step size (which multiples the covariance in the proposal) at the start of a run to achieve proposal acceptance rates near to 25\%. (This is near optimal in a precise sense for certain high-dimensional posteriors \cite{rob}; 
in practice, it works well beyond this setting.) All sampling is performed in the transformed space where the prior distribution is normal. Figure~\ref{fig:mcmc} shows kernel density estimates of the MCMC results; the panels in the left column are for EKI-GP (for $\myvec{y}^{(1)},\dots,\myvec{y}^{(4)}$), and the panels in the right column are for B-GP for the same data realization. We display contours of the posterior that contain $50\%,75\%,$ and $99\%$ of the mass of the posterior density.

All sets of results converge to similar regions of the parameter space about the true parameters, and the spread of uncertainty is quantified similarly in both EKI-GP and B-GP.  Table \ref{table:stddev} shows the standard deviations of the individual parameters alongside the empirical standard deviation calculated from the ensemble spread in EKI iteration 9. The standard deviations from the MCMC posterior based on EKI-GP and B-GP are similar. {We re-emphasize
that the EKI is constructed as an ensemble optimization method and has the property that the ensemble 
evolves towards consensus among the parameters, while
also matching the data: the ensemble collapses. As a result, the EKI ensemble spread is an inadequate estimate the uncertainty, as 
seen in Table \ref{table:stddev}.  
As explained in the introduction, ensemble methods are only justifiable to quantify uncertainties in the Gaussian posterior setting \cite{oli,rey}.}
Our approach is justifiable whenever
the GP accurately approximates the
forward model \cite{Cleary21a}. The use of
EKI for the design of training points for the
GP does not require accurate uncertainty quantification
within EKI; it only relies on EKI approximately locating
minimizers of the model-data misfit.

There is sampling variability because of the different data realizations. %(One sees similar perturbations of the red + representing the ensemble mean of the final EKI optimization iteration.)
This sampling variability can be assessed by asking 
which probability contours contain the true parameters. For both EKI-GP and B-GP, in three of four realizations, we capture the true values within 50\% of the posterior probability mass; the realization $\myvec{y}^{(3)}$ is captured only within the 99\% contour of the posterior probability.

\subsection{Uncertainty Quantification in Prediction Experiments} \label{sec:forwardUQ}
\begin{figure}
 
   \centering
   \includegraphics[width=0.7\textwidth]{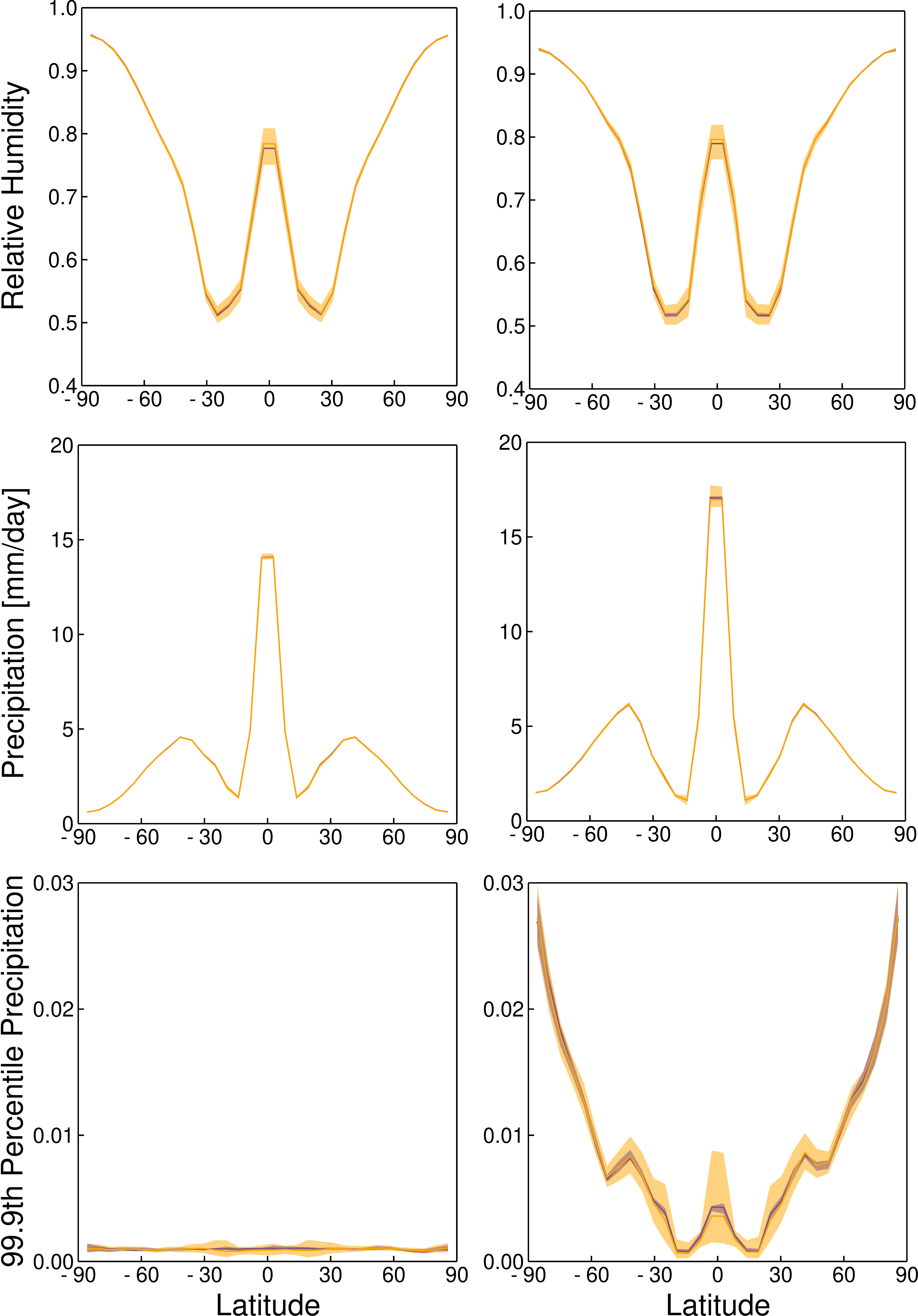}
   
  \caption{Comparison of statistics of a $7200$-day average in a climate-change simulation. Left column: control climate; right column: warmer climate. Synthetic observational data evaluated at the true fixed parameters are shown in blue, while data evaluated at 100 samples from the posterior distribution (EKI-GP) are shown in orange. (We choose the posterior from the first data realization, top-left panel of Figure \ref{fig:mcmc}.) The solid lines are the medians, and the shaded regions represent the 95\% confidence intervals (between the 2.5th and 97.5th percentiles). Top: Relative humidity in mid-troposphere. Middle: Precipitation rate. Bottom: Frequency with which 99.9th percentile of latitude-dependent daily precipitation in the control climate is exceeded.}
  \label{fig:forward_uq_all}
\end{figure}
%Validation of the posterior distribution using emulated forward model requires running MCMC directly with the GCM, over very long time-horizons; this is not feasible with our computationally expensive forward model. \hl{[obscure--people will not understand this]}
To illustrate how the posterior distribution of parameters obtained in the sample step of the CES algorithm can be used to produce climate predictions with quantified uncertainties, we consider an idealized global-warming experiment. As in \cite{OGorman08c,OGorman08b}, we rescale the longwave opacity of the atmosphere everywhere by a uniform factor $\alpha$. In the control climate we have considered so far, $\alpha = 1 $. We generate a warm climate by setting $\alpha=1.5$, which results in a global-mean surface air temperature increase from 287~K in the control climate to 294~K in the warm climate. To see parametric uncertainty rather than internal variability noise in the resulting ``climate change predictions,'' we use long (7,200-day or  approximately 20-year) averages in the prediction experiments. 

We evaluate predictions of the latitude-dependent relative humidity and mean precipitation rate that we used in the CES algorithm. We also consider the frequency of precipitation extremes, now taken as the frequency with which the 99.9th percentile of daily precipitation in the control simulation is exceeded (rather than the 90th percentile we considered earlier). This last statistic indicates how the frequency of what are 1-in-1000 day precipitation events in the control climate change in the warmer climate. 

We investigate the effect of parametric uncertainty on predictions by taking 100 samples of parameters from the posterior, creating a prediction for each sample, and comparing statistics of these runs with runs in which parameters are fixed to the true values $\myvec{\theta}^\dagger$. The climate statistics in the control climate are shown in the left column of Figure \ref{fig:forward_uq_all}. The runs from posterior samples (orange) and with fixed true parameters (blue) match well. The noise due to internal variability is quantitatively represented by the blue shaded region. Unlike in the earlier figures with short (30-day) averages (e.g., Figure \ref{fig:validation_benchmark}), the internal variability noise here is small relative to the parametric uncertainty because of the (long) 7200-day averaging window. The orange shaded region contains both internal variability and parametric uncertainty and is dominated by parametric uncertainty. This remains the case in the warmer climate (right column of Figure \ref{fig:forward_uq_all}). 

The effects of global warming on atmospheric quantities is seen by comparing the two columns of Figure \ref{fig:forward_uq_all}. Relative humidity is fairly robust to the warming climate, and precipitation rates increase globally \cite{OGorman08b}. The most dramatic changes occur for the frequency of extreme precipitation events \cite{OGorman09a}. What is a 1-in-1000 day event in the control climate (e.g., occuring with frequency 0.001) occurs in the extratropics of the warmer climate an order of magnitude more frequently, with the 95\% confidence interval spanning $0.01$ to $0.03$. That is, a 1-in-1000 day event in the control climate occurs  every $30$ to $100$ days in the warmer climate. The parametric uncertainty is particularly large for extreme precipitation events within the tropics---behavior one would not be able to see in global warming experiments with fixed parameters. This is consistent with the known high uncertainty in predictions of tropical rainfall extremes with comprehensive climate models \cite{OGorman09b}.

\section{Conclusion and Discussion} \label{sec:discussion}

The primary goal of this article was to demonstrate that ensemble Kalman inversion (EKI), machine
learning, and MCMC algorithms can be judiciously combined
within the calibrate-emulate-sample framework to efficiently estimate uncertainty of model parameters in computationally expensive climate models. We provided a proof-of-concept in a relatively simple idealized GCM. 

Our approach is novel because we train a machine learning (GP) emulator using input-output pairs generated from an EKI algorithm. This methodology has several advantageous features:
\begin{enumerate}
\item It requires a modest number of runs of the expensive forward model (typically, $O(100)$ runs). 
\item It generally finds optimal or nearly optimal parameters even in the presence of internal variability noise because EKI is robust with respect to such noise. 
\item The resulting GP emulation is naturally most accurate around the (a priori unknown) optimal parameters because this is where EKI training points concentrate. 
\item MCMC shows robust convergence to the posterior distribution, and allows identification of the optimal parameters with the maximum of the posterior probability, because it uses an objective function that is smoothed by GP emulation. 
\end{enumerate}
The effectiveness of GP depends on the training points, and a user must choose how many iterations of EKI to use for training (before ensemble collapse). In practice, we find the GP performance is robust as long as we include the initial iteration of training points (drawn from the prior) in our training set. The necessity of using the initial
ensemble could be side-stepped by using an ensemble method that does not collapse, such as the recently introduced ensemble Kalman sampler (EKS) \cite{GarHofLiStu20}.

The CES algorithm is efficient, as it addresses two dominant sources of computational expense. First, poor prior knowledge of model parameters requires blind exploration of a possibly high-dimensional parameter space
to find optimal parameters and thus the region of high posterior probability. The CES framework handles this with an EKI algorithm, which we show to be successful when using time averaged data from a chaotic nonlinear model. Second, computing parametric uncertainty with a sampling technique (such as MCMC) generally requires many ($10^5$--$10^6$) evaluations of an expensive forward model. We instead solve a cheap approximate problem by exploiting GP emulators. We train the emulators on relatively few ($O(100)$) intelligently chosen evaluations provided by EKI, which ensures that training points are placed where they are most needed---near the minimum of the model-data misfit. The training itself introduces negligible computational cost relative to the running of the forward model, and the computational expense of evaluating the emulator in the sampling step is also negligible. Hence, the CES framework achieves about a factor 1000 speedup over brute-force MCMC algorithms. Significant efforts to accelerate brute-force MCMC without approximation have been undertaken \cite{Jarvinen10a,Solonen2012}, and improvements of up to a factor 5 speedup have been made with adaptive and parallelized Markov chains. However, these approaches still are considerably more expensive than the CES algorithm.

The CES algorithm also has a smoothing property, which is beneficial even in situations where a forward model is cheap enough to apply a brute-force MCMC. If the forward model exhibits internal variability, the objective function for the sampling algorithm will contain a data misfit of the form \eqref{eq:objfunc}, which has a random component because it contains a finite-time average. Without more sophisticated sampling methods, MCMC algorithms get stuck in noise-induced local minima. In the CES algorithm, only EKI uses the functional \eqref{eq:objfunc}, and EKI is well suited for this purpose. The GP emulator learns the smooth, noiseless  model $\mathcal{G}_\infty$ (in which internal variability disappears), using evaluations of $\mathcal{G}_T$ (which are affected by internal variability). Thus, MCMC within the CES algorithm uses the smooth GP approximation of  \eqref{eq:objfunc_inf}.

The MCMC results in this study successfully capture the true
parameters and their uncertainties. The results contain natural biases arising from the use of prior distributions, internal variability of the climate, and use of a
single noisy sample as synthetic data. Despite the sampling variability and emulator constraints, our
MCMC samples were able to capture the true parameters within an estimated 99\% confidence
interval in our examples, demonstrating the potential of
EKI-trained GP emulators for MCMC sampling. Validation of the emulator (Figure \ref{fig:validation}) further supports the MCMC results,
as do our comparisons with MCMC using the benchmark emulator
(Table \ref{table:stddev}). The GP emulator both smooths the objective function and allows us to quantify uncertainty by sampling from the posterior distribution.  This contrasts with uncertainty quantification based on the EKI ensemble, which underestimates the true uncertainties in our experiments by an order of magnitude. As used here, EKI should be viewed as an optimization algorithm and not a sampling algorithm.
Adding additional spread to match the posterior
within EKI may be achieved for
Gaussian posteriors \cite{oli,rey} or by means of EKS \cite{GarHofLiStu20}; however, these methods
are not justifiable beyond the Gaussian setting. The
MCMC algorithm within CES, on the other hand, samples from an
approximate posterior distribution and is justifiable
beyond the Gaussian posterior setting \cite{Cleary21a}.

\revtwo{The key to the success of the CES methodology is founded
on ensemble Kalman methods and the cheap derivative-free
manner in which they can learn states or, in the setting
of this paper, parameters, from data. We emphasize that
there exists already significant work showing, empirically,
that ensemble Kalman methods can be generalized to perform
uncertainty quantification beyond the Gaussian
setting \cite{bocquet2012combining,gu2007iterative,li2007iterative,sakov2012,bocquet2014,defforge2019improving,defforge2021improving}
by means of iterated EnKF, which may be motivated via the Newton
method \cite{sakov2012}, and also by means of ensemble-based Gauss-Newton
\cite{sakov2018iterative} in the presence of model error.
However, it is unclear how this methodology might be rigorously justified and developed as a general-purpose
tool, despite the empirical successes it has shown on a number of
challenging problems. In contrast,
CES may be systematically refined, based on the theory of
Gaussian process emulation \cite{williams2006gaussian} and its use
within Bayesian inference \cite{stuart2018posterior}. CES does not depend on the
UQ properties of ensemble methods, only on their optimization
and estimation properties, which are amenable to systematic analysis
\cite{Schillings17a,kelly2014well,de2018long}. For these reasons, CES has considerable potential as a systematic approach to uncertainty quantification.}

Good scaling of the CES algorithm with the dimension of the parameter space will be of critical importance for moving beyond the current, low-dimensional proof-of-concept setting. Each stage of the CES algorithm scales to higher dimensions: For the calibration stage, ensemble methods scale well to high-dimensional state and parameter spaces, typically with $O(10^2)$ forward model runs \cite{Kalnay03a,Oliver08a}, if used with localization. In high dimensions, regularization is also typically needed in calibration algorithms, and various regularization schemes can be added  \cite{chada2020tikhonov,iglesias2015iterative,iglesias2016regularizing,Garbuno-Inigo2019,schneider2020imposing}.
The sampling stage also scales well to high dimensions. In particular, MCMC methods scale to non-parametric (infinite-dimensional) problems, 
in which the unknown is a function \cite{Cotter13a}. 
The current bottleneck for scalability lies with the GP emulator, which does not scale easily to high-dimensional inputs. However, other supervised machine learning techniques have potential to do so. For example, random feature maps and deep neural networks show promise in this regard \cite{nelsen2020random,bhattacharya2020model}; incorporating these tools in the CES algorithm is a direction of current research. 

 An alternative form of constraining parameter uncertainty is history matching, or precalibration \cite{BowGolVer10,CamEdwRou11,Wil_etal13}. The idea complements that of Bayesian uncertainty quantification, where instead of searching for a high probability region of parameter space with respect to data, one rules out regions of parameter space that are deemed inconsistent with the data.  \cite{Cou_etal20_pre} and \cite{Hou_etal20_pre} recently constrained the parameter space of a parameterization scheme by approximating a plausibility function over the parameter space using a Gaussian process, and then removing ``implausible" regions of parameter space where the plausibility function passes a threshold. This removal process is iterated until the uncertainty of the emulator is small enough, or the space becomes empty. History matching accomplishes a similar adaptivity task as performed in the CES algorithm by EKI. During early stages of history matching, however, one must sample the full parameter space with reasonable resolution, and emulator training is required at every iteration to evaluate the plausibility function. In contrast, in the CES algorithm, EKI draws a modest numbers of samples at every iteration and can work directly with noisy model evaluations, lowering the computational expense. The output of history matching is a (possibly empty) acceptable set of forward model runs; sampling this set leads to an upper bound on the prediction uncertainty. The benefit of the CES algorithm is that it provides samples of the posterior distribution, which lead to full estimates of prediction uncertainty (Figure  \ref{fig:forward_uq_all}). For this reason, history matching has been proposed as a preprocessing step for Bayesian uncertainty quantification, known as precalibration  to improve priors and assess model validity \cite{BowGolVer10,CamEdwRou11}. The CES algorithm targets the Bayesian posterior distributions directly. %, where reduction of parameter space without the need for observational data, and the possibility of empty solution space \cite{BowGolVer10} could be used to . 

In the more comprehensive climate modeling settings we target, data will be given from earth observations and from local high-resolution simulations \cite{Schneider17c}. In these settings, model error leads to deficiencies when comparing model evaluations to data, leading to structural biases and uncertainty that must be quantified. Structural model errors can be quantified with a flexible hierarchical Gaussian process regression that learns a non-parametric form of the model deficiency, as demonstrated in prototype problems in \cite{jin2}. %We represent the model-data residual as a stochastic differential equation (SDE). Using the statistics of finite time averaged data, a GP is used to learn a non-parametric form of the deficiency by modeling drift and diffusion of this SDE. 
This approach represents model error in an interpretable
fashion, as part of the model itself, rather than in the data space as pioneered
in \cite{Kennedy01a}.

The CES framework has potential for both the calibration and uncertainty quantification of comprehensive climate models, and other computationally expensive models. It is computationally efficient enough to use data averaged in time (e.g., over seasons), which need to be accumulated over longer model runs. Time-averaged climate statistics, including mean values and higher-order statistics such as extreme value statistics, are what typically matters in climate predictions. CES allows us to focus model calibration and uncertainty quantification on such immediately relevant statistics. Using time averaged statistics also has the advantage that it leads to smoother, albeit still noisy, objective functions when compared with calibration of climate models by minimizing mismatches in instantaneous, short-term forecasts \cite{Schneider17c}. The latter approach can improve short-term forecasts but may not translate into improved climate simulations \cite{Schirber2013}. It also suffers from the difficulty that model resolution and data resolution may be mismatched. Focusing on climate statistics, as we did in our proof-of-concept here, circumvents this problem: time-aggregated climate statistics are varying relatively smoothly in space and, hence, minimizing mismatches in statistics between models and data does not suffer from the resolution-mismatch problem. CES can be used to learn about arbitrary parameters in climate models from time-averaged data. It leads to quantification of parametric uncertainties that then can be converted into parametric uncertainties in predictions by sampling from the posterior distribution of parameters.

\bigskip\noindent
\textbf{Acknowledgements.}
This work was supported by the generosity of Eric and Wendy Schmidt by recommendation of the Schmidt Futures program, by the Hopewell Fund, the Paul G. Allen Family Foundation, and the National Science Foundation (NSF, award AGS‐1835860). A.M.S. was also supported by the Office of Naval Research (award N00014-17-1-2079). We thank Emmet Cleary for his preliminary work underlying some of the results shown here. \revtwo{We would like to thank the reviewers for their insightful comments and suggestions which have lead to the improvement of this article.}

\bigskip\noindent
\textbf{Data Availability.} All computer code used in this paper is open source. The code for the idealized GCM, the Julia code for the CES algorithm, the plot tools, and the slurm/bash scripts to run both GCM and CES are available at https://doi.org/10.5281/zenodo.4393029. 
%%  --  --  --  --  --  --  --  --  --  --  --  --  --  --  --  --  --  --  --  --  --  --  --  --  %%
%% References and Citations

%%%%%%%%%%%%%%%%%%%%%%%%%%%%%%%%%%%%%%%%%%%%%%%
%
% \bibliography{<name of your .bib file>} don't specify the file extension
%
% don't specify bibliographystyle
%%%%%%%%%%%%%%%%%%%%%%%%%%%%%%%%%%%%%%%%%%%%%%%
\bibliographystyle{siam}
\bibliography{add,lit,library}

\end{document}